\documentclass[11pt,epsfig]{article}
\usepackage{amsmath}
\usepackage{amsthm}
 \usepackage{amssymb}
\usepackage{epsfig}
\usepackage{graphicx}
\usepackage{amsfonts}
\usepackage{color, epstopdf}
\usepackage{cite}
\usepackage{graphicx,wrapfig}
\usepackage{flafter}
\usepackage{fancyhdr}
\usepackage{stmaryrd}
\usepackage{graphicx,subfigure}
\usepackage{multicol}
\usepackage{dsfont}
\usepackage{booktabs}
\usepackage[center]{caption2}
\usepackage{epstopdf}
\usepackage [latin1]{inputenc}
\usepackage{multirow}
\usepackage{enumerate}
\usepackage{enumitem}
\usepackage{algpseudocode,algorithm,algorithmicx}

\newtheorem{theorem}{Theorem}[section]
\newtheorem{lemma}{Lemma}[section]

\newtheorem{example}{Example}

\newcommand{\diff}{\triangledown_{\!\tau}}

\newcommand{\defeq}{:=}
\newcommand{\zd}{\,\mathrm{d}}
\newcommand{\abs}[1]{\left|#1\right|}
\newcommand{\absb}[1]{\big|#1\big|}

\newcommand{\bra}[1]{\left(#1\right)}
\newcommand{\brab}[1]{\big(#1\big)}
\newcommand{\braB}[1]{\Big(#1\Big)}
\newcommand{\brat}[1]{(#1)}
\newcommand{\kbra}[1]{\left[#1\right]}
\newcommand{\kbrab}[1]{\big[#1\big]}
\newcommand{\kbraB}[1]{\Big[#1\Big]}

\newcommand{\mynormb}[1]{\big\|#1\big\|}

\title{On energy stable, maximum-principle preserving, second order BDF scheme with variable steps for the Allen-Cahn equation}
\author{Hong-lin Liao\thanks{
Department of Mathematics,
Nanjing University of Aeronautics and Astronautics,
Nanjing 211106, P. R. China. E-mails: liaohl@nuaa.edu.cn and liaohl@csrc.ac.cn.
This author's work is supported by a grant 1008-56SYAH18037 from NUAA Scientific Research Starting Fund of Introduced Talent.}
\quad Tao Tang\thanks{SUSTech International Center for Mathematics, Shenzhen, China; and
    Division of Science and Technology, BNU-HKBU United International College,
    Zhuhai, Guangdong Province, China.
    Email: tangt@sustech.edu.cn. This author's work is partially supported by the NSF of China under grant number 11731006 and the science challenge project (No. TZ2018001),}
\quad Tao Zhou\thanks{NCMIS \& LSEC, Institute of Computational Mathematics and Scientific/Engineering Computing,
Academy of Mathematics and Systems Science, Chinese Academy of Sciences, Beijing, 100190,
P. R. China. Email: tzhou@lsec.cc.ac.cn. This author's work is partially supported by the NSF of China (under grant numbers 11822111, 11688101, 11571351, and 11731006), the science challenge project (No. TZ2018001), NCMIS, and the youth innovation promotion association (CAS).}
}
\date{}

\begin{document}

\maketitle

\begin{abstract}
In this work, we investigate the two-step backward differentiation formula (BDF2) with nonuniform grids for the Allen-Cahn equation. We show that the nonuniform BDF2 scheme is energy stable under the time-step ratio restriction $r_k:=\tau_k/\tau_{k-1}<(3+\sqrt{17})/2\approx3.561.$ Moreover, by developing a novel kernel recombination and complementary technique, we show, for the first time, the discrete maximum principle of BDF2 scheme under the time-step ratio restriction $r_k<1+\sqrt{2}\approx 2.414$ and a practical time step constraint.  The second-order rate of convergence in the maximum norm is also presented. Numerical experiments are provided to support the theoretical findings.
\\[1ex]

\end{abstract}

\section{Introduction}
\setcounter{equation}{0}
The phase field equations are important models in describing a host of free-boundary problems in various areas,
including material, physical and biology systems \cite{Allen1979A,Cahn1958Free,WiseLFC:2008,HVO:2006}. Meanwhile, numerical schemes for phase field equations have also been extensively studied in recent years \cite{Gomez2011Provably,ChenWangWangWise:2014,HeLiuTang:2006,MaChenYangZhang:2017,YanChenWangWise:2018,WangWangWise:2010}. The main focuses of the numerical schemes are the discrete energy stability (e.g., \cite{WangWangWise:2010,ChenWangWangWise:2014,ChenWangYanZhang:2019,ShenYang:2010,ShenYangYu:2015}) and the discrete maximum principle (for Allen-Cahn equations) \cite{Hou2017Numerical,TangYuZhou:2019,Du2019} which are inherent properties in the continuous level. Another key feature of the phase field models is that the associate solutions in general admit multiple time scales, i.e. an initial dynamics evolves on a fast time scale and later coarsening evolves on a very slow time scale. This motives the use of nonuniform meshes in time domain \cite{Gomez2011Provably,LuoTangXie:2016,QiaoAn2011,Zhang2013An,Li2017Computationally,FengTangYang:2015}, i.e., one adopts small time steps to capturing the fast dynamics when the solution varies rapidly
while uses large time steps otherwise to accelerate the time integration. While the numerical analysis for numerical schemes with uniform grids has been well investigated, however, the relevant analysis for nonuniform grids have not been well studied. In fact, even for linear/semilinear parabolic equations, the relevant study is far from complete \cite{Becker:1998,Emmrich:2005,LeRoux:1982}.

To this end, we investigate in this work the well known two-step backward differentiation formula (BDF2) \cite{Becker:1998,Emmrich:2005,GearTu:1974,Grigorieff:1983,LeRoux:1982,Nishikawa:2019,YanChenWangWise:2018} with nonuniform grid for the Allen-Cahn equation. As a simple phase field model, the Allen-Cahn equation admits the energy dissipation law and the maximum principle in the continuous level, and our purpose is to investigate whether the nonuniform BDF2 scheme can preserve these properties in the discrete level. Compared to existing literature, our contributions are three folds:
\begin{itemize}
\item We show that the nonuniform BDF2 scheme is energy stable under the time-step ratio restriction
$r_k:=\tau_k/\tau_{k-1}<(3+\sqrt{17})/2\approx3.561.$
\item We show, for the first time, the discrete maximum principle of nonuniform BDF2 scheme under
the time-step ratio restriction $r_k<1+\sqrt{2}$ and a practical time step constraint.
\item We show the second-order rate of convergence in the maximum norm, and present several experiments to support the theoretical findings.
\end{itemize}
We mention a related work \cite{ChenWangYanZhang:2019}, where the nonuniform BDF2 scheme (combined with the convex splitting approach) is investigated for the Cahn--Hilliard equation, and the energy stability and convergence analysis are presented under similar time-step ratio restrictions as in the current work. The key tool in \cite{ChenWangYanZhang:2019} for the optimal error estimates is a generalized discrete Gronwall inequality. In contrast, we develop in this work a novel kernels recombination and complement (KRC) technique for the analysis. Moreover, our proof for the discrete maximum principle of nonuniform BDF2 scheme seems to be the first work with such results.

The rest of this paper is organized as following. In Section 2, we provide with some preliminaries. The discrete maximum principle and the discrete energy stability are presented in Section 3 and Section 4, respectively. In Section 5, we show the rigorous convergence analysis in the  maximum norm, and this is followed by several numerical examples in Section 6. We finally give some concluding remarks in Section 7.

\section{Preliminaries}
We consider the following Allen-Cahn equation:
\begin{align}
\partial_{t}u(\mathbf{x},t)
=&\,\varepsilon^{2}\Delta{u}-f(u),\quad
\mathbf{x}\in\Omega, \quad 0<{t}\leq{T},\label{Problem-1}\\
u(\mathbf{x},0)=&\,u_{0}(\mathbf{x}),\quad \mathbf{x}\in\bar{\Omega},\label{Problem-2}
\end{align}
where $\mathbf{x}=(x,y)^{T}$ and $\Omega=(0,L)^2$ with its closure $\bar{\Omega}$.
The nonlinear bulk force $f(u)$ is given by $f(u)=u^{3}-u$, and the small constant  $0< \varepsilon \ll 1 $ is
the interaction length that describes
the thickness of the transition boundary between materials. For simplicity, we consider the periodic boundary conditions.
As is well known, the above Allen-Cahn equation can be viewed an $L^2$-gradient flow of
the following Ginzburg--Landau free energy functional
\begin{align}\label{eq: continuous energy}
E[u](t):=\int_{\Omega}\braB{\frac{1}{2}\varepsilon^{2}|\nabla{u}|^{2}
+F[u]}\zd{\mathbf{x}}, \quad F[u]=\frac{1}{4}(1-u^{2})^{2}.
\end{align}
In other words, the Allen-Cahn equation \eqref{Problem-1} admits the following energy dissipation law
\begin{align}\label{eq: continuous energy law}
\frac{\zd E}{\zd t}\leq0.
\end{align}
Moreover, the following maximum principle holds
\begin{align}\label{eq: continuous Maximum Principle}
|u(\mathbf{x},t)|\leq{1}, \quad
\text{if}\quad|u(\mathbf{x},0)|\leq{1}.
\end{align}
\subsection{The nonuniform BDF2 scheme}
We consider a general nonuniform time grid $0=t_0<t_1<t_2<\cdots<t_N=T$ with the time-step
$\tau_k\defeq t_k-t_{k-1}$ for $1\le k\le N$, and the maximum step size
$\tau\defeq\max_{1\le k\le N}\tau_k$.
For any time sequence $\{v^n\}_{n=0}^N$, we denote $\diff v^n:=v^n-v^{n-1}$ and $\partial_{\tau}v^n:=\diff v^n/\tau_n$.
For $k=1,2$, let $\Pi_{n,k}v$ be the interpolating polynomial of a function~$v$
over $k+1$ nodes $t_{n-k}$, $\cdots$, $t_{n-1}$~and $t_{n}$.
Then by taking $v^n=v(t_n)$, the BDF1 formula yields
$$D_1v^n:=\bra{\Pi_{n,1}v}'(t)=\diff  v^n/\tau_n, \quad n\ge1,$$
and furthermore, the well known BDF2 formula reads
\begin{align}\label{eq: BDF2 formula}
D_2v^n:=\bra{\Pi_{n,2}v}'(t_n)=&\,\frac{1+2r_{n}}{\tau_n(1+r_{n})}\diff v^{n}
-\frac{r_{n}^2}{\tau_n(1+r_{n})}\diff  v^{n-1},\quad  n\ge2,
\end{align}
where the adjacent time-step ratios $r_k$ are defined by $r_1\equiv0$ (if necessary) and
$$r_{k}\defeq\frac{\tau_{k}}{\tau_{k-1}}, \quad  2\le k\le N.$$
To introduce the fully discrete scheme, we consider a central finite difference approximation in physical domain. For a positive integer $M$,
let $h:=L/M$ be the spatial grid length
and we set $\bar{\Omega}_{h}:=\big\{\mathbf{x}_{h}=(ih,jh)\,|\,0\leq i,j\leq M\}$.
For any grid function $\{v_h\,|\,\mathbf{x}_{h}\in\bar{\Omega}_{h}\}$, we denote
$$\mathbb{V}_{h}:=\big\{v\,|\,v=(v_{j})^{T}\;\;\text{for}\;\;1\leq{j}\leq{M},
\;\;\text{with}\;\;
v_{j}=(v_{i,j})^{T}\;\;\text{for}\;\;1\leq{i}\leq{M}\big\},$$
where $v^{T}$ is the transpose of the vector $v$.
We also define the associate maximum norm $\|v\|_{\infty}:=\max_{\mathbf{x}_{h}\in\Omega_{h}}|v_{h}|$.
We shall denote by $\Lambda_h$ the discrete matrix of Laplace operator $\Delta$ subject to periodic boundary conditions.

In general, one can use the BDF1 scheme to obtain first-level solution $u^1$ by considering $D_2v^1:=D_1v^1,$ as the two-step BDF2 formula needs two starting values and the BDF1 scheme generates a second-order accurate solution at the first time grid. Then, we have the following fully discrete nonlinear BDF2 time-stepping scheme
\begin{align}\label{eq: BDF2 vector scheme}
D_2u^n
&=\varepsilon^{2}\Lambda_hu^{n}-f(u^{n}), \quad   n\geq{1},
\end{align}
where the vector $f(u^{n})$ is defined in the element-wise, that is, $f(u^{n}):=(u^{n})^{.3}-u^{n}$.
\subsection{Summary of main contributions}
The main purpose of this work is to analyze the nonuniform BDF2 scheme \eqref{eq: BDF2 vector scheme}. In particular, we shall show in Theorem \ref{thm: enery stability} in the next section that scheme \eqref{eq: BDF2 vector scheme} admits
a discrete energy stability,
under the following mild time-step ratio constraint
\begin{enumerate}[itemindent=1em]
  \item[\textbf{S1}.] $0< r_k< \brat{3+\sqrt{17}}/2\approx3.561, \,\,\,2\le k\le N$.
\end{enumerate}
Then, we present the discrete maximum principle and convergence estimates of scheme \eqref{eq: BDF2 vector scheme} in Section 4. To do this, we shall propose a novel kernels recombination and complementary (KRC) technique. More precisely, the BDF2 formula
\eqref{eq: BDF2 formula} is first regarded as a discrete convolution summation,
\begin{align}\label{eq: BDF2}
D_2v^n=\sum_{k=1}^nb^{(n)}_{n-k}\diff  v^k, \quad  n\ge1,
\end{align}
where the discrete convolution kernels $b^{(n)}_{n-k}$ are defined by $b^{(1)}_0:=1/\tau_1$ and
\begin{align}\label{eq: BDF2 kernels}
b^{(n)}_0:=\frac{1+2r_n}{\tau_n(1+r_n)},\quad b^{(n)}_1:=-\frac{r_n^2}{\tau_n(1+r_n)}
\quad\text{and}\quad b^{(n)}_j:=0\;\; \textmd{for} \;\; 2\leq j\leq n.
\end{align}
For notation simplicity, we set $b^{(n)}_n:=0$ for $n\ge1$ when necessary, and set $\sum_{k=i}^{j} \cdot =0 $ if the index $i>j$.

In the \emph{kernels recombination stage} of KRC, we introduce a new class of variables $\{\bar{v}^k\}$
that consist of a linear combination of the original variables $\{v^k\}$
and reformulate $D_2v^n$ into a new discrete convolution form, such as $\sum_{k=1}^nd^{(n)}_{n-k}\diff \bar{v}^k$
which always involves all of previous solutions $\{\bar{v}^{k}\}_{k=0}^{n-1}$.
The main aim is to build a new class of discrete convolution kernels $d^{(n)}_{j}$ so that
they are nonnegative and monotonously decreasing.
Then we show in Theorem \ref{thm:Dis-Max-Principle} in Section 4 that
the scheme \eqref{eq: BDF2 vector scheme} preserves the maximum principle
under a time-step ratio restriction that coincides with the zero-stability condition due to Grigorieff \cite{Grigorieff:1983}:
\begin{enumerate}[itemindent=1em]
  \item[\textbf{S0}.] $0< r_k<1+\sqrt{2}\approx2.414$ for $2\le k\le N$.
\end{enumerate}

The discrete maximum principle offers us the possibility to show the maximum norm convergence without any Lipschitz assumptions on the nonlinear bulk force.
With the help of the \emph{kernels complementary stage} of KRC,
we build in Lemma \ref{lem: gronwall} a new discrete Gr\"{o}nwall inequality. Then we show in Theorem \ref{thm: maxmimum norm convergence} that
the scheme \eqref{eq: BDF2 vector scheme} is of second-order rate of convergence in the maximum norm under the step-ratio condition \textbf{S0}.
To the best of our knowledge, it is the \emph{first} work establishing such
convergence results for nonuniform BDF2 scheme under the Grigorieff's zero-stability condition \textbf{S0}.

\section{Solvability and energy stability}
\setcounter{equation}{0}
We first list some well known properties of the matrix $\Lambda_h$ in the following lemma.
\begin{lemma}\label{lem:Negative-Condition}
The discrete matrix $\Lambda_h$ of Laplace operator $\Delta$ has the following properties
\begin{itemize}
  \item [(a)] The discrete matrix $\Lambda_h$ is symmetric.
  \item [(b)] For any nonzero $v\in{\mathbb{V}_{h}}$, $v^{T}\Lambda_hv\leq{0}$, i.e., the matrix $\Lambda_h$ is negative semi-definite.
  \item [(c)] The elements of $\Lambda_h=(d_{ij})$ fulfill $d_{ii}=-\max_{i}\sum_{j\neq{i}}|d_{ij}|$ for each $i$.
\end{itemize}
\end{lemma}

Then, we show the solvability of scheme \eqref{eq: BDF2 vector scheme} in the following lemma.
\begin{lemma}\label{lem:solvability}
The discrete scheme \eqref{eq: BDF2 vector scheme} is uniquely solvable
if  $$\tau_n<\frac{1+2r_{n}}{1+r_{n}}, \quad n \geq 1.$$
Notice that the above step constraint is practical as it is suffice to require $\tau_n< 1.$
\end{lemma}
\begin{proof}We rewrite the nonlinear scheme \eqref{eq: BDF2 vector scheme} into
\begin{align*}
G_hu^{n}+(u^{n})^{.3}=g(u^{n-1})\quad \textmd{with}\quad  g(u^{n-1}):=b_0^{(n)}u^{n-1}-b_{1}^{(n)}\diff u^{n-1}, \quad n\ge1,
\end{align*}
where $G_h:=b_0^{(n)}-1-\varepsilon^{2}\Lambda_h$.
If the time-step size $\tau_n<\frac{1+2r_{n}}{1+r_{n}}$, by definition \eqref{eq: BDF2 kernels} we have $b_0^{(n)}>1$.
Thus the matrix $G_h$ is positive definite according to Lemma \ref{lem:Negative-Condition} (b).
Consequently, the solution of nonlinear equations solves
  \begin{align*}
u^n=\arg \min_{w\in \mathbb{V}_{h}}\left\{\frac12w^TG_hw+\frac{1}4\sum_{k=1}^{M}w_k^{4}-w^Tg(u^{n-1})\right\}, \quad n\geq1.
\end{align*}
The strict convexity of the above objective function implies the unique solvability of \eqref{eq: BDF2 vector scheme}.
\end{proof}

We now consider the energy stability of the nonuniform BDF2 scheme \eqref{eq: BDF2 vector scheme} by
defining a modified discrete energy $\widehat{E}:$
\begin{align}\label{eq: discrete energy}
\widehat{E}[u^k]:=E[u^k]+\frac{r_{k+1}\tau_k}{2(1+r_{k+1})}\sum_{i=1}^{M}\brab{\partial_\tau u_i^k}^2, \quad k\ge1,
\end{align}
where we set $\widehat{E}[u^0]=E[u^0]$ which corresponds to the setting $r_1\equiv0$, and $E[u^k]$ is the original discrete energy that is given by
\begin{align*}
E[u^k]:=-\frac{\varepsilon^{2}}2\brat{u^{k}}^T\Lambda_hu^{k}+\frac14\sum_{i=1}^{M}\brab{1-(u_i^k)^{2}}^2, \quad k\ge0.
\end{align*}
Notice that the modified energy $\widehat{E}[u^k] \rightarrow E[u^k]$ when $\tau \rightarrow 0.$ We are now ready to present the following energy stability of scheme \eqref{eq: BDF2 vector scheme}.
\begin{theorem}\label{thm: enery stability}
Assume that the step-ratio condition \textbf{S1} holds, and moreover, suppose that
\begin{align}\label{condition: discrete energy law}
\tau_k\le \min\bigg\{\frac{1+2r_{k}}{1+r_{k}},\frac{2+4r_k-r_k^2}{1+r_k}-\frac{r_{k+1}}{1+r_{k+1}}\bigg\}\quad\text{for $k\ge1$.}
\end{align}
Then, the discrete solution $u^n$ of the BDF2 time-stepping scheme \eqref{eq: BDF2 vector scheme} satisfies
\begin{align}\label{eq: discrete energy law}
	\widehat{E}[u^k]\le \widehat{E}[u^{k-1}], \quad k\ge1.
\end{align}
\end{theorem}
\begin{proof}
Taking the $L^2$ inner product (in the vector space) of \eqref{eq: BDF2 vector scheme} with $\brat{\diff u^n}^T$, we have
\begin{align}\label{proof: energy equality}
\sum_{i=1}^{M}D_2u_i^{n}\brab{\diff u_i^n}
	-\varepsilon^2\brat{\diff u^{n}}^T\Lambda_h u^{n}
        +\sum_{i=1}^{M}f(u_i^n)\diff u_i^n=0,\quad  n\ge1.
\end{align}
By using  Lemma \ref{lem:Negative-Condition} (a)-(b), one gets
\begin{align*}
-\varepsilon^{2}\brat{\diff u^n}^T\Lambda_hu^{n}=&\,
-\frac{\varepsilon^{2}}2\brat{u^{n}}^T\Lambda_hu^{n}+\frac{\varepsilon^{2}}2\brat{u^{n-1}}^T\Lambda_hu^{n-1}
-\frac{\varepsilon^{2}}2\brat{\diff u^n}^T\Lambda_h\brab{\diff u^n}\\
\ge&\,-\frac{\varepsilon^{2}}2\brat{u^{n}}^T\Lambda_hu^{n}+\frac{\varepsilon^{2}}2\brat{u^{n-1}}^T\Lambda_hu^{n-1}.
\end{align*}
It is easy to check the following identity
\begin{align*}
4\bra{a^3-a}\bra{a-b}+2(1-a^2)\bra{a-b}^2=&\,(1-a^2)^2-(1-b^2)^2+(a^2-b^2)^2.
\end{align*}
Taking $a:=u_i^n$ and $b:=u_i^{n-1}$ in the above equality we obtain
\begin{align*}
\sum_{i=1}^{M}f(u_i^{n})\brab{\diff u_i^n}
=&\,\sum_{i=1}^{M}\brab{(u_i^n)^{3}-u_i^{n}}\brab{\diff u_i^n}\\
\ge&\,\frac14\sum_{i=1}^{M}\brab{1-(u_i^n)^{2}}^2
-\frac14\sum_{i=1}^{M}\brab{1-(u_i^{n-1})^{2}}^2
-\frac{1}{2}\sum_{i=1}^{M}\brab{\diff u_i^n}^2.
\end{align*}
Thus it follows from \eqref{proof: energy equality} that
\begin{align}\label{proof: energy equality2}
\sum_{i=1}^{M}D_2u_i^{n}\brab{\diff u_i^n}-\frac{\tau_n^2}{2}\sum_{i=1}^{M}\brab{\partial_{\tau} u_i^n}^2
	+E(u^n)\le E(u^{n-1}),\quad   n\ge1.
\end{align}
We now consider the mathematical induction argument. For the case of $n=1$, we have
\begin{align*}
D_2u_i^{1}\brab{\diff u_i^1}=&\,D_1u_i^{1}\brab{\diff u_i^1}=\frac{r_{2}\tau_1}{2(1+r_{2})}\brab{\partial_\tau u_i^1}^2
+\frac{2+r_2}{2(1+r_{2})}\tau_1\brab{\partial_\tau u_i^1}^2\\
\ge&\,\frac{r_{2}\tau_1}{2(1+r_{2})}\brab{\partial_\tau u_i^1}^2
+\frac{\tau_1^2}{2}\brab{\partial_\tau u_i^1}^2,
\end{align*}
where the condition \eqref{condition: discrete energy law} of $k=1$ was used in the last inequality.
The estimate \eqref{proof: energy equality2} then gives
$$\widehat{E}[u^1]\le \widehat{E}[u^{0}]=E[u^{0}].$$
For the general case of $n\ge2$, we use the identity $2a(a-b)=a^2-b^2+(a-b)^2$ and the definition \eqref{eq: BDF2 kernels} of BDF2 kernels to obtain
\begin{align*}
D_2u_i^{n}\brab{\diff u_i^n}&=\brab{b_0^{(n)}+b_1^{(n)}}\brab{\diff u_i^n}^2-b_1^{(n)}\bra{\diff u_i^n-\diff u_i^{n-1}}\diff u_i^n\\
&=\brab{b_0^{(n)}+\frac12b_1^{(n)}}\brab{\diff u_i^n}^2+\frac12b_1^{(n)}\brab{\diff u_i^{n-1}}^2
-\frac12b_1^{(n)}\bra{\diff u_i^n-\diff u_i^{n-1}}^2\\
&\ge\brab{b_0^{(n)}+\frac12b_1^{(n)}}\brab{\diff u_i^n}^2+\frac12b_1^{(n)}\brab{\diff u_i^{n-1}}^2\\
&=\frac{r_{n+1}\tau_n}{2(1+r_{n+1})}\brab{\partial_{\tau}u_i^n}^2
-\frac{r_{n}\tau_{n-1}}{2(1+r_{n})}\brab{\partial_{\tau}u_i^{n-1}}^2
+\braB{\frac{2+4r_n-r_n^2}{1+r_n}-\frac{r_{n+1}}{1+r_{n+1}}}\frac{\tau_n}{2}\brab{\partial_{\tau}u_i^n}^2.
\end{align*}
Inserting this estimate into \eqref{proof: energy equality2}, we obtain
\begin{align*}
\braB{\frac{2+4r_n-r_n^2}{1+r_n}-\frac{r_{n+1}}{1+r_{n+1}}-\tau_n}\frac{\tau_n}{2}\sum_{i=1}^{M}\brab{\partial_{\tau}u_i^n}^2
+\widehat{E}[u^n]\le \widehat{E}[u^{n-1}],\quad  2\le n\le N.
\end{align*}
The desired result follows by noticing the restriction \eqref{condition: discrete energy law}, and this completes the proof.
\end{proof}

Some comments for the time-step restriction \eqref{condition: discrete energy law} are listed below.
The first constraint in \eqref{condition: discrete energy law} comes from Lemma \ref{lem:solvability} for solvability, and one is suffice to choose $\tau_k\le1$ to ensure it for any $r_k>0$.

It remains to check the second constraint in \eqref{condition: discrete energy law}.
For $n=1$, the constraint \eqref{condition: discrete energy law} yields $\tau_1\le \frac{2+r_2}{1+r_2}$ and
one can also simply choose $\tau_1\le1$. Under the condition \textbf{S1}, one has $0<r_k< r_s$, where $r_s=\frac{3+\sqrt{17}}{2}$
is the positive root of the algebraic equation $2+3r_s-r_s^2=0$,
and $\frac{r_{k+1}}{1+r_{k+1}}<\frac{r_s}{1+r_s}=\frac{\sqrt{17}-1}{4}\approx0.78$.
So the time-step restriction \eqref{condition: discrete energy law} are fulfilled by choosing
\begin{align*}
\tau_k\le\frac{2+4r_k-r_k^2}{1+r_k}-\frac{r_{s}}{1+r_{s}}=\frac{2+4r_k-r_k^2}{1+r_k}-\frac{\sqrt{17}-1}{4}\quad\text{for $k\ge2$.}
\end{align*}
Actually, let $h(x):=\frac{2+4x-x^2}{1+x}$ such that
$h'(x)=\frac{x+1+\sqrt{3}}{(1+x)^{2}}\brat{\sqrt{3}-1-x}.$
We consider three cases:
\begin{itemize}
  \item[(i)] If $0<r_k\le \sqrt{3}-1$, then $h'(r_k)\ge0$ and $h(r_k)\ge h(0)=2$. One can choose time-steps $\tau_k\le\min\big\{1,\frac{9-\sqrt{17}}{4}\big\}=1$ to ensure \eqref{condition: discrete energy law}.
  \item[(ii)] If $\sqrt{3}-1<r_k\le \sqrt{2}+1$, then $h'(r_k)<0$ and $h(r_k)\ge h(\sqrt{2}+1)=1+\frac{\sqrt{2}}2$. One can choose time-steps $\tau_k\le1+\frac{\sqrt{2}}2-\frac{\sqrt{17}-1}{4}\approx0.93$ to ensure \eqref{condition: discrete energy law}.
  \item[(iii)] If $\sqrt{2}+1<r_k< r_s$, then $h'(r_k)< 0$ and $h(r_k)> h(r_s)=\frac{r_s}{1+r_s}$.
  In this case, especially when the current step-ratio $r_k\rightarrow r_s$, one can choose a small time-step
  $\tau_{k+1}$ or step-ratio $r_{k+1}$ to ensure the time-step restriction \eqref{condition: discrete energy law}
  in adaptive computations. For an example, the time-steps $\tau_k\le\frac12$ are sufficient if one choose
  the next step-ratio $r_{k+1}\le \frac{2h(r_s)-1}{3-2h(r_s)}\approx0.39$.
\end{itemize}

To summary, under the condition \textbf{S1}, the time-step size constraint \eqref{condition: discrete energy law} is reasonable. In particular, it is practical in controlling the next time-step $\tau_{k+1}$ in adaptive simulations.

\section{Kernels recombination and discrete maximum principle}
\setcounter{equation}{0}
In this section, we shall show the discrete maximum principle of scheme \eqref{eq: BDF2 vector scheme}.
\subsection{Reformation of BDF2 formula}
We first introduce a new class of variables below (see \cite[Remark 6]{LiaoMcLeanZhang:2019} for technical motivations):
\begin{align}\label{eq: BDF2 newvariable}
\bar{v}^0:=v^0\quad\text{and}\quad \bar{v}^k:=v^k-\eta v^{k-1}\;\;\text{for $k\geq1$},
\end{align}
where $\eta$ is a real parameter to be determined.
It is easy to find the substitution formula
\begin{align}\label{eq: substitution formulas}
v^k=\bar{v}^k+\eta v^{k-1}=\bar{v}^k+\eta\brab{\bar{v}^{k-1}+\eta v^{k-2}}=\cdots
=\sum_{\ell=0}^{k}\eta^{k-\ell}\bar{v}^{\ell}\quad
\text{for $k\geq1,$}
\end{align}
and then we have
\begin{align*}
\diff v^k=\sum_{\ell=1}^{k}\eta^{k-\ell}\diff\bar{v}^{\ell}+\eta^{k}v^0\quad
\text{for $k\geq1.$}
\end{align*}
By inserting the above equation into \eqref{eq: BDF2} and exchanging the summation order, we obtain an updated BDF2 formula
\begin{align}\label{eq: BDF2 new form}
D_2v^n\equiv\sum_{j=1}^nd^{(n)}_{n-j}\diff\bar{v}^j+d^{(n)}_{n}\bar{v}^0\quad\text{for $n\ge1$},
\end{align}
where the new discrete convolution kernels $d^{(n)}_{n-j}$ can be defined by
\begin{align*}
d^{(n)}_{n-j}:=\sum_{k=j}^{n}b^{(n)}_{n-k}\eta^{k-j}\quad\text{for $1\leq j\leq n$},
\quad\text{and}\quad d^{(n)}_{n}:=\eta d^{(n)}_{n-1}.
\end{align*}
Alternatively, we have the following explicit formula
\begin{align}\label{eq: BDF2 new kernels}
d^{(n)}_{0}:=b_{0}^{(n)}\quad\text{and}\quad d^{(n)}_{j}:=\eta^{j-1}\brab{b_{0}^{(n)}\eta+b_{1}^{(n)}}
\quad\text{for $1\leq j\leq n$}.
\end{align}
We shall require that the new discrete kernels $d^{(n)}_{n-j}$ are nonnegative and decreasing, that is, $d^{(n)}_{0}\ge d^{(n)}_{1}\ge\cdots\ge d^{(n)}_{n}\ge0$. By the definitions \eqref{eq: BDF2 new kernels} and \eqref{eq: BDF2 kernels}, it is easy to check that
this aim can be achieved by setting
\begin{align}\label{eq: ETA condition}
\frac{r_k^2}{1+2r_k}\le \eta< 1\quad\text{for $k\ge2$}.
\end{align}
Meanwhile, we require that the adjacent time-step ratios satisfy the condition \textbf{S0}, that is, $r_k<1+\sqrt{2},$
which coincides with the Grigorieff's zero-stability condition \cite{Grigorieff:1983} for ODE problems.

Now, by using the new formula \eqref{eq: BDF2 new form}, the numerical scheme \eqref{eq: BDF2 vector scheme} reads
\begin{align}\label{eq: alternaltive BDF2 scheme0}
\sum_{j=1}^nd^{(n)}_{n-j}\diff\bar{u}^j+d^{(n)}_{n}\bar{u}^0
&=\varepsilon^{2}\Lambda_hu^{n}-f(u^{n})\quad\text{for $n\geq{1}$}.
\end{align}
This equation will be our starting point to establish the discrete maximum principle.
Recalling the definition of $\bar{u}^j$ and the substitution formula
\eqref{eq: substitution formulas}, we have
\begin{align}\label{eq: alternaltive BDF2 scheme1}
\brab{d^{(n)}_{0}-1-\varepsilon^{2}\Lambda_h}u^{n}+(u^{n})^{.3}
=&\,\eta d^{(n)}_{0}u^{n-1}+\sum_{j=0}^{n-1}\brab{d^{(n)}_{n-j-1}-d^{(n)}_{n-j}}\bar{u}^j\nonumber\\
=&\,d^{(n)}_{0}\sum_{j=0}^{n-1}\eta^{n-j}\bar{u}^{j}+\sum_{j=0}^{n-1}\brab{d^{(n)}_{n-j-1}-d^{(n)}_{n-j}}\bar{u}^j\quad\text{for $n\geq{1}$}.
\end{align}
This formulation \eqref{eq: alternaltive BDF2 scheme1} will be used to evaluate $u^n$ by using the information from $\big\{\bar{u}^j\big\}_{j=0}^{n-1}$. Again, we apply the substitution formula
\eqref{eq: substitution formulas} to derive from  \eqref{eq: alternaltive BDF2 scheme0}
that
\begin{align}\label{eq: alternaltive BDF2 scheme2}
\brab{d^{(n)}_{0}+S_n-\varepsilon^{2}\Lambda_h}\bar{u}^n
=&\,\sum_{j=0}^{n-1}\brab{d^{(n)}_{n-j-1}-d^{(n)}_{n-j}-S_n\eta^{n-j}+\eta^{n-j}\varepsilon^{2}\Lambda_h}\bar{u}^j\nonumber\\
&\,+(S_n+1)u^{n}-(u^{n})^{.3}\nonumber\\
=&\,\sum_{j=0}^{n-1}Q_{n-j}^{(n)}\bar{u}^j+(S_n+1)u^{n}-(u^{n})^{.3}\quad\text{for $n\geq{1}$},
\end{align}
where $S_n$ is a real parameter (that can depend on the time-levels) to be determined, and the matrix
\begin{align}\label{eq: matrix Qj}
Q_j^{(n)}:=\brab{d^{(n)}_{j-1}-d^{(n)}_{j}-S_n\eta^{j}}I+\eta^{j}\varepsilon^{2}\Lambda_h\quad\text{for $1\le j\le n$.}
\end{align}
This formulation will be used to evaluate $\bar{u}^n$ by using the information $\big\{\bar{u}^j\big\}_{j=0}^{n-1}$ and $u^n$.

\subsection{Choice of recombined parameter}

Next lemma presents a time-step size restriction so that
the matrix $Q_j^{(n)}$ in \eqref{eq: matrix Qj} is bounded in the maximum norm.

\begin{lemma}\label{lem: Matrix-Qj-Inf-Norm}
Assume that the step-ratio condition \textbf{S0} holds, and suppose that the time-step size satisfies
\begin{align}\label{eq: taun condition}
\tau_n\le \frac{(1+2r_n)\eta-r_n^2}{\eta^2(1+r_n)}\frac{1-\eta}{S_n+4\varepsilon^{2}h^{-2}}\quad
\text{for $n\ge1$},
\end{align}
where the recombined parameter $\eta$ satisfies \eqref{eq: ETA condition}.
Then the matrix $Q_j^{(n)}$ in \eqref{eq: matrix Qj} fulfills
\begin{align}\label{eq: Qj estimate}
\mynormb{Q_j^{(n)}}_{\infty}\le d^{(n)}_{j-1}-d^{(n)}_{j}-S_n\eta^{j}\quad\text{for $1\le j\le n$.}
\end{align}
\end{lemma}
\begin{proof} Consider the case of $n\ge2$. By the definition \eqref{eq: BDF2 new kernels}, the matrix $Q_j^{(n)}$ in \eqref{eq: matrix Qj} reads
\begin{align*}
Q_j^{(n)}=\eta^{j}\kbra{(1-\eta)\eta^{-2}\brab{b_{0}^{(n)}\eta+b_{1}^{(n)}}-S_n}I+\eta^{j}\varepsilon^{2}\Lambda_h\quad\text{for $2\le j\le n$.}
\end{align*}
The time-step condition \eqref{eq: taun condition} together with the definition \eqref{eq: BDF2 kernels} yields
\begin{align*}
\frac{1-\eta}{\eta^2}\brab{b_{0}^{(n)}\eta+b_{1}^{(n)}}- S_n\ge\frac{4\varepsilon^{2}}{h^{2}}\,.
\end{align*}
Thus all the elements of the matrix $Q_j^{(n)}=\brab{q_{k\ell}^{(n,j)}}$ are nonnegative and
\begin{align*}
\mynormb{Q_{j}^{(n)}}_{\infty}=\max_{k}\sum_{\ell}\absb{q_{k\ell}^{(n,j)}}
=\max_{k}\sum_{\ell}q_{k\ell}^{(n,j)}\le d^{(n)}_{j-1}-d^{(n)}_{j}-S_n\eta^{j}\quad\text{for $2\le j\le n$.}
\end{align*}
The desired estimate \eqref{eq: Qj estimate} holds for $2\le j\le n$.
It remains to consider the case $j=1$ for $n\ge1$.
By using the step condition \eqref{eq: taun condition}, the definitions \eqref{eq: BDF2 kernels}
and \eqref{eq: BDF2 new kernels} show that (with $r_1=0$)
\begin{align*}
d^{(n)}_{0}-d^{(n)}_{1}-S_n\eta=&\,(1-\eta)b_{0}^{(n)}-b_{1}^{(n)}-S_n\eta\\
=&\,\eta^{-1}\kbrab{(1-\eta)\brab{b_{0}^{(n)}\eta+b_{1}^{(n)}}-b_{1}^{(n)}-S_n\eta^2}\\
\ge&\,\eta\kbraB{(1-\eta)\eta^{-2}\brab{b_{0}^{(n)}\eta+b_{1}^{(n)}}-S_n}\ge \frac{4\eta\varepsilon^{2}}{h^{2}}.
\end{align*}
Thus, all elements of the matrix $Q_1^{(n)}=\brab{q_{k\ell}^{(n,1)}}$ are nonnegative and
\begin{align*}
\mynormb{Q_{1}^{(n)}}_{\infty}=\max_{k}\sum_{\ell}\absb{q_{k\ell}^{(n,1)}}
=\max_{k}\sum_{\ell}q_{k\ell}^{(n,1)}\le d^{(n)}_{0}-d^{(n)}_{1}-S_n\eta.
\end{align*}
The proof is complete.
\end{proof}
Further comments for the restriction \eqref{eq: taun condition} are listed below. We set
$$K(\eta):=\frac{1-\eta}{\eta^2}\frac{(1+2r_n)\eta-r_n^2}{1+r_n}.$$
Obviously, $K(\eta)>0$ if
the parameter $\eta$ satisfies \eqref{eq: ETA condition}. Moreover,
$K'(\eta)=\frac{1+r_n}{\eta^3}\brab{\frac{2r_n^2}{(1+r_n)^2}-\eta},$
and $K(\eta)$ approaches its maximum value when $\eta\rightarrow\frac{2r_n^2}{(1+r_n)^2}.$
For a fixed maximum step-ratio $r_s\in[1, 1+\sqrt{2})$,
one can choose the parameter $\eta\in\big[\frac{r_s^2}{1+2r_s},1)$ such that
the condition \eqref{eq: ETA condition} holds at any time-levels.
To relieves the restriction \eqref{eq: taun condition} on the time-step size,
we can choose in all above derivations
\begin{align}\label{eq: ETA Choice}
\eta:=\frac{2r_s^2}{(1+r_s)^2} \quad \text{with $r_s\in[1, 1+\sqrt{2})$}
\end{align}
For example, consider the uniform mesh case with $r_n=r_s=1$, one can take $\eta=\frac12$ so that
the time-step condition \eqref{eq: taun condition} reads
$$
\tau_n=\tau\le \frac{1}{2(S_n+4\varepsilon^{2}h^{-2})}.$$
Consider the case of $r_s=2$, one can take the recombined parameter $\eta=8/9$
so that the time-step condition \eqref{eq: taun condition} requires
$$\tau_n\le \frac{1}{48}\frac{1}{S_n+4\varepsilon^{2}h^{-2}}.$$
The time-step condition \eqref{eq: taun condition} with $S_2=2$ will be used to establish
the discrete maximum principle in next subsection.

\subsection{Discrete maximum principle}

To establish the discrete maximum-principle, we recall the following result \cite[Lemma 3.2]{Hou2017Numerical}.
\begin{lemma}\label{lem:Matrix-Inf-Norm}
Let $B$ be a real $M\times{M}$ matrix and $A=aI-B$ with $a>0$.
If the elements of $B=(b_{ij})$ fulfill $b_{ii}=-\max_{i}\sum_{j\neq{i}}|b_{ij}|$,
then for any $c>0$ and $V\in{\mathbb{R}^{M}}$ we have
\begin{align}
\|AV\|_{\infty}\geq{a}\|V\|_{\infty}\quad\text{and}\quad
\|AV+c(V)^{3}\|_{\infty}\geq{a}\|V\|_{\infty}+c\|V\|_{\infty}^{3}.\nonumber
\end{align}
\end{lemma}
We are now ready to present the following theorem on discrete maximum principle.
\begin{theorem}\label{thm:Dis-Max-Principle}
Assume that the step-ratio restriction \textbf{S0} holds and suppose that the time-step size satisfies
\begin{align}\label{eq: tau condition-Max-Principle}
\tau_n\le \frac{(1+2r_n)\eta-r_n^2}{\eta^2(1+r_n)}\frac{1-\eta}{2+4\varepsilon^{2}h^{-2}}\quad
\text{for $n\ge1$},
\end{align}
where the recombined parameter $\eta$ is defined by \eqref{eq: ETA Choice}.
Then, the BDF2 time-stepping scheme \eqref{eq: BDF2 vector scheme}
preserves the maximum principle at the discrete levels,
that is,
\begin{align*}
\mynormb{u^{k}}_{\infty}\leq 1\;\;\text{for $1\leq{k}\leq{N}$}\quad\text{if\; $\mynormb{u^{0}}_{\infty}\leq{1}$.}
\end{align*}
\end{theorem}
\begin{proof}
The desired result is a by-product of the following claim
\begin{align*}
\mynormb{\bar{u}^{k}}_{\infty}\leq 1-\eta\;\;\text{for $1\leq{k}\leq{N}$}\quad\text{if\; $\mynormb{\bar{u}^{0}}_{\infty}\leq{1}$.}
\end{align*}
We now verify this new claim with the complete mathematical induction argument.
Taking $n=1$ in \eqref{eq: alternaltive BDF2 scheme1}, one has
\begin{align*}
\brab{d^{(1)}_{0}-1-\varepsilon^{2}\Lambda_h}u^{1}+(u^{1})^{.3}
=&\,\eta d^{(1)}_{0}u^{0}+(1-\eta)d^{(1)}_{0}\bar{u}^0=d^{(1)}_{0}\bar{u}^0.
\end{align*}
Since $d^{(1)}_{0}=b_0^{(n)}>1$, we apply Lemmas \ref{lem:Negative-Condition} and \ref{lem:Matrix-Inf-Norm} to get
\begin{align*}
\brab{d^{(1)}_{0}-1}\mynormb{u^{1}}_{\infty}+\mynormb{u^{1}}_{\infty}^{3}
\le&\,\mynormb{\brab{d^{(1)}_{0}-1-\varepsilon^{2}\Lambda_h}u^{1}+(u^{1})^{.3}}_{\infty}\le d^{(1)}_{0},
\end{align*}
which implies $\mynormb{u^{1}}_{\infty}\leq1$.
To see this, notice that the function $g_c(z):=\bra{c-1}z+z^{3}-c$ is increasing with respect to $z>0$, if the real parameter $c\ge1$.
So this contradicts with  $\mynormb{u^{1}}_{\infty}>1.$

Next we shall bound $\mynormb{\bar{u}^{1}}_{\infty}$.
Because $|(c+1)z-z^3|\le c$ for $\abs{z}\le1$
if  the real parameter $c\ge2$, one has
$\mynormb{3u^{1}-(u^{1})^{.3}}_{\infty}\le 2$.
Thus we take $n=1$ and $S_1=2$ in the equation \eqref{eq: alternaltive BDF2 scheme2} and
apply Lemma \ref{lem: Matrix-Qj-Inf-Norm} to get
\begin{align*}
\brab{d^{(1)}_{0}+2}\mynormb{\bar{u}^1}_{\infty}\le &\,
\mynormb{\brab{d^{(1)}_{0}+2-\varepsilon^{2}\Lambda_h}\bar{u}^1}_{\infty}
=\mynormb{Q_{1}^{(1)}\bar{u}^{0}+3u^{1}-(u^{1})^{.3}}_{\infty}\\
\le&\,\mynormb{Q_{1}^{(1)}}_{\infty}\mynormb{\bar{u}^{0}}_{\infty}+\mynormb{3u^{1}-(u^{1})^{.3}}_{\infty}\\
\le &\,d^{(1)}_{0}-d^{(1)}_{1}-2\eta+2=(1-\eta)\brab{d^{(1)}_{0}+2},
\end{align*}
which yields $\mynormb{\bar{u}^{1}}_{\infty}\leq1-\eta$.

For the general case of $2\leq n\le N$, assume that
\begin{align}\label{inductionAssume}
\mynormb{\bar{u}^{k}}_{\infty}\leq1-\eta\quad\text{for $1\leq{k}\leq{n-1}.$}
\end{align}
From the equation \eqref{eq: alternaltive BDF2 scheme1} and the expressions in \eqref{eq: BDF2 new kernels}, one applies
Lemmas \ref{lem:Negative-Condition} and \ref{lem:Matrix-Inf-Norm} to find
\begin{align*}
\brab{d^{(n)}_{0}-1}\mynormb{u^{n}}_{\infty}+\mynormb{u^{n}}_{\infty}^{3}
\le&\,\mynormb{\brab{d^{(n)}_{0}-1-\varepsilon^{2}\Lambda_h}u^{n}+(u^{n})^{.3}}_{\infty}\nonumber\\
\le&\,d^{(n)}_{0}\sum_{j=0}^{n-1}\eta^{n-j}\mynormb{\bar{u}^{j}}_{\infty}
+\sum_{j=0}^{n-1}\brab{d^{(n)}_{n-j-1}-d^{(n)}_{n-j}}\mynormb{\bar{u}^{j}}_{\infty}\\
\le&\,\eta d^{(n)}_{0}
+(1-\eta)\brab{d^{(n)}_{0}-d^{(n)}_{n-1}}+\brab{d^{(n)}_{n-1}-d^{(n)}_{n}}=d^{(n)}_{0},
\end{align*}
where the inductive hypothesis \eqref{inductionAssume} and the identity
$(1-\eta)\sum_{j=1}^{n-1}\eta^{n-j}+\eta^{n}=\eta$ have been used in the third inequality.
This yields immediately
\begin{align}\label{eq: Un estimate}
\mynormb{u^{n}}_{\infty}\leq1.
\end{align}
It remains to evaluate $\mynormb{\bar{u}^{n}}_{\infty}$.
The above estimate \eqref{eq: Un estimate} gives 
\begin{align*}
\mynormb{3u^{n}-(u^n)^{.3}}_{\infty}\le 2.
\end{align*}
Now we take $S_n=2$ in the equation \eqref{eq: alternaltive BDF2 scheme2}.
By applying Lemma \ref{lem: Matrix-Qj-Inf-Norm} and
the inductive hypothesis \eqref{inductionAssume} one has
\begin{align*}
\brab{d^{(n)}_{0}+2}\mynormb{\bar{u}^n}_{\infty}
\le&\,\sum_{j=0}^{n-1}\mynormb{Q_{n-j}^{(n)}}_{\infty}\mynormb{\bar{u}^j}_{\infty}+\mynormb{3u^{n}-(u^n)^{.3}}_{\infty}\\
\le&\,(1-\eta)
\sum_{j=1}^{n-1}\brab{d^{(n)}_{n-j-1}-d^{(n)}_{n-j}-2\eta^{n-j}}
+\brab{d^{(n)}_{n-1}-d^{(n)}_{n}-2\eta^{n}}+2\\
=&\,(1-\eta)\brab{d^{(n)}_{0}-d^{(n)}_{n-1}}+\brab{d^{(n)}_{n-1}-d^{(n)}_{n}}
-2(1-\eta)\sum_{j=1}^{n-1}\eta^{n-j}-2\eta^{n}+2\\
=&\,(1-\eta)\brab{d^{(n)}_{0}+2}.
\end{align*}
This leads to $\mynormb{\bar{u}^n}_{\infty}\le 1-\eta$, and the proof is completed.
\end{proof}

Notice that in the Allen-Cahn equation \eqref{Problem-1},
the coefficient $\varepsilon\ll 1$ represents the width of diffusive interface. In practice,
one should choose a small spacial step $h=\mathcal{O}(\varepsilon)$ to track the moving interface.
Then the restriction \eqref{eq: tau condition-Max-Principle} is approximately equivalent to
\begin{align*}
\tau_n\le \frac{(1+2r_n)\eta-r_n^2}{\eta^2(1+r_n)}\frac{1-\eta}{6}\quad
\text{for $\eta:=\frac{2r_s^2}{(1+r_s)^2}$ \,\, and \,\, $n\ge1$}.
\end{align*}
On the other hand,
the parameter $\eta$ is introduced only for
the theoretical analysis but not necessary in numerical computations,
thus the time-step restriction \eqref{eq: tau condition-Max-Principle} seems to be rather practical. We also remark that Theorem 4.1 seems to be the first result on second order maximum-principle preserving scheme with variable steps.

\section{Complementary kernels and convergence analysis}
\setcounter{equation}{0}
This section is devoted to convergence analysis. To this end, we introduce a class of discrete complementary convolution kernels $\big\{(Q_{\!d})_{n-j}^{(n)}\big\}_{j=1}^n$
via the discrete kernels $d_{j}^{(n)}$ in \eqref{eq: BDF2 new kernels},
\begin{align}\label{eq: convolutionKernels def}
(Q_{\!d})_{0}^{(n)}:=\frac{1}{d_{0}^{(n)}}\quad\text{and}\quad
(Q_{\!d})_{n-j}^{(n)}:=
\sum_{k=j+1}^{n}\frac{d_{k-j-1}^{(k)}-d_{k-j}^{(k)}}{d_0^{(j)}}(Q_{\!d})_{n-k}^{(n)}
    \quad \text{for $1\leq j\leq n-1$}.
\end{align}
This type of discrete kernels was first introduced in \cite{LiaoLiZhang:2018}
for numerical approximation of fractional Caputo derivatives and
further generalized in \cite{LiaoMcLeanZhang:2019} for more general discrete kernels.
It is easy to check that the following complementary identity holds
\begin{align}
\sum_{j=k}^n(Q_{\!d})^{(n)}_{n-j}d_{j-k}^{(j)}\equiv 1
\quad\text{for $\forall\;1\le k\le n$.}\label{eq: complementary identity}
\end{align}
From the definition \eqref{eq: BDF2 new kernels}, we know that $d_{j}^{(n)}$ are nonnegative and decreasing.
So the definition \eqref{eq: convolutionKernels def} implies that $(Q_{\!d})_{n-j}^{(n)}\ge0$. The identity
\eqref{eq: complementary identity} yields immediately
\begin{align}
0<(Q_{\!d})^{(n)}_{n-j}\le \frac{1}{d^{(j)}_{0}}
\quad\text{for\;\; $\forall\;1\le j\le n$.}\label{eq: complementary estimate}
\end{align}

Now we apply the discrete complementary convolution kernels $\big\{(Q_{\!d})_{n-j}^{(n)}\big\}_{j=1}^n$
and their properties \eqref{eq: complementary identity}-\eqref{eq: complementary estimate}
to build a novel discrete Gr\"{o}nwall lemma, which will plays an important role for the analysis of the nonuniform BDF2 scheme.
\begin{lemma}\label{lem: gronwall}
For constants $\kappa>0$, $\lambda\in(0,1)$ and for any non-negative sequences~$\{g^k\}_{k=1}^N$ and $\{w^k\}_{k=0}^N$ such that
\begin{align*}
\sum_{k=1}^nd^{(n)}_{n-k}\diff w^k\le \kappa\sum_{k=1}^n\lambda^{n-k}w^{k}+g^n
    \quad\text{for $1\le n\le N$,}
\end{align*}
where the discrete kernels $d^{(n)}_{j}$ are defined by \eqref{eq: BDF2 new kernels}.
If $b^{(n)}_{0}\ge2\kappa$, then
\begin{align*}
w^n\le2\exp\brab{\frac{2\kappa t_n}{1-\lambda}}\braB{w^0+\sum_{j=1}^n\frac{g^j}{b^{(j)}_{0}}}\quad\text{for $1\le n\le N$.}
\end{align*}
\end{lemma}
\begin{proof}We have
\begin{align*}
\sum_{k=1}^jd^{(j)}_{j-k}\diff w^k\le \kappa\sum_{k=1}^j\lambda^{j-k}w^{k}+g^j
    \quad\text{for $1\le j\le N$.}
\end{align*}
Multiplying the above inequality by the complementary kernels
$(Q_{\!d})^{(n)}_{n-j}$ and taking the index $j$ from $1$ to $n$ one gets
\begin{align*}
\sum_{j=1}^n(Q_{\!d})^{(n)}_{n-j}\sum_{k=1}^jd^{(j)}_{j-k}\diff w^k\le \kappa\sum_{j=1}^n(Q_{\!d})^{(n)}_{n-j}\sum_{k=1}^j\lambda^{j-k}w^{k}+\sum_{j=1}^n(Q_{\!d})^{(n)}_{n-j}g^j.
\end{align*}
By exchanging the summation order and applying the complementary identity
\eqref{eq: complementary identity}, one has
\begin{align*}
&\sum_{j=1}^n(Q_{\!d})^{(n)}_{n-j}\sum_{k=1}^jd^{(j)}_{j-k}\diff w^k=
\sum_{k=1}^n\diff w^k\sum_{j=k}^n(Q_{\!d})^{(n)}_{n-j}d^{(j)}_{j-k}=w^n-w^0,\\
&\sum_{j=1}^n(Q_{\!d})^{(n)}_{n-j}\sum_{k=1}^j\lambda^{j-k}w^{k}
=\sum_{k=1}^nw^{k}\sum_{j=k}^n(Q_{\!d})^{(n)}_{n-j}\lambda^{j-k}.
\end{align*}
Thus it follows that
\begin{align*}
w^n\le&\, w^0+\kappa w^{n}(Q_{\!d})^{(n)}_{0}+
2\kappa\sum_{k=1}^{n-1}w^{k}\sum_{j=k}^n(Q_{\!d})^{(n)}_{n-j}\lambda^{j-k}+\sum_{j=1}^n(Q_{\!d})^{(n)}_{n-j}g^j
\quad\text{for $1\le n\le N$.}
\end{align*}
Furthermore, the estimate \eqref{eq: complementary estimate} and the definition \eqref{eq: BDF2 new kernels} yields
$$(Q_{\!d})^{(n)}_{n-1}\leq \frac1{b^{(1)}_{0}}=\tau_1\quad\text{and}\quad
(Q_{\!d})^{(n)}_{n-j}\leq \frac1{b^{(j)}_{0}}=\frac{1+r_j}{1+2r_j}\tau_j\le \tau_j\quad \text{for $2\leq j\leq n$}.$$
Setting $b^{(n)}_{0}\ge2\kappa$ so that $(Q_{\!d})^{(n)}_{0}\le\frac1{b^{(n)}_{0}}\le\frac1{2\kappa}$, then one gets
\begin{align*}
w^n\le&\, 2w^0+
2\kappa\sum_{k=1}^{n-1}w^{k}\sum_{j=k}^n(Q_{\!d})^{(n)}_{n-j}\lambda^{j-k}+2\sum_{j=1}^n(Q_{\!d})^{(n)}_{n-j}g^j\\
\le&\, 2\kappa\sum_{k=1}^{n-1}w^{k}\sum_{j=k}^n\tau_j\lambda^{j-k}
+2w^0+2\sum_{j=1}^n\frac{g^j}{b^{(j)}_{0}}
\quad\text{for $1\le n\le N$.}
\end{align*}
Note that
$$4\sum_{k=1}^{n-1}\sum_{j=k}^n\tau_j\lambda^{j-k}\le 4\sum_{j=1}^{n}\tau_j\sum_{k=1}^j\lambda^{j-k}\le \frac{4t_n}{1-\lambda}.$$
The desired result follows by the standard Gr\"{o}nwall inequality and the proof is completed.
\end{proof}
We are now ready to present the following convergence result:
\begin{theorem}\label{thm: maxmimum norm convergence}
Let the initial data $u_0$ be smooth and bounded by $1$,
and the solution of \eqref{Problem-1}-\eqref{Problem-2} be sufficiently smooth.
Assume that the step-ratio restriction \textbf{S0} holds and the time-step size satisfies
\eqref{eq: tau condition-Max-Principle}.
The numerical solution $u_{h}^{n}$
of the BDF2 scheme \eqref{eq: BDF2 vector scheme} is convergent in the maximum norm, and it holds
\begin{align*}
\mynormb{u(\mathbf{x}_h,t_n)-u_{h}^{n}}_{\infty}\leq \frac{C_ut_n}{1-\eta}\exp\left({\frac{4t_n}{1-\eta}}\right)\bra{\tau^2+h^{2}}\quad\text{for $1\leq{n}\leq{N},$}
\end{align*}
where the recombined parameter $\eta$ is determined by \eqref{eq: ETA Choice}, and $C_u$ is a constant that is independent of the time-step sizes and time-step ratios.
\end{theorem}
\begin{proof} Let $U_{h}^{n}:=u(\mathbf{x}_h,t_n)$ and $e_{h}^{n}:=U_{h}^{n}-u_{h}^{n}\in{\mathbb{V}_{h}}$
for $\mathbf{x}_{h}\in\bar{\Omega}_{h}$ and $0\leq{n}\leq{N}$.
It is easy to find that the exact solution $U_{h}^{n}$ satisfies the governing equation
\begin{align*}
D_2U^{n}=\varepsilon^{2}\Lambda_{h}U^{n}-f(U^{n})+\Upsilon^{n}+R^{n},\quad {1}\leq{n}\leq{N},
\end{align*}
where $\Upsilon^{n}$ and $R^{n}$ denote the truncation errors in time and space, respectively.
Subtracting the numerical scheme \eqref{eq: BDF2 vector scheme} from the above equation one gets
\begin{align}\label{CONproof: ErrorEquation}
D_2e^{n}=\varepsilon^{2}\Lambda_{h}e^{n}+f(u^{n})-f(U^{n})+\Upsilon^{n}+R^{n},\quad{1}\leq{n}\leq{N}
\end{align}
with $e^0=0$. As done before, we define $\bar{e}^k:=e^k-\eta e^{k-1}$ for $k\ge1$ with $\bar{e}^0:=e^0=0$.
Recalling the elementary inequality
$$|(a^{3}-a)-(b^{3}-b)|\leq{2}|a-b|\quad\text{for $\forall\, a,b\in[-1,1]$},$$
we apply Theorem \ref{thm:Dis-Max-Principle} (discrete maximum principle) to get
\begin{align}\label{CONproof: nonlinear term}
\mynormb{f(U^{n})-f(u^{n})}_{\infty}\leq 2\mynormb{e^{n}}_{\infty}.
\end{align}
By using the alternative formulas \eqref{eq: substitution formulas}-\eqref{eq: BDF2 new form},
we rewrite the error equation \eqref{CONproof: ErrorEquation} into
\begin{align*}
\sum_{j=1}^nd^{(n)}_{n-j}\diff\bar{e}^j-\varepsilon^{2}\sum_{j=1}^n\eta^{n-j}\Lambda_{h}\bar{e}^j
=f(u^{n})-f(U^{n})+\Upsilon^{n}+R^{n},\quad{1}\leq{n}\leq{N},
\end{align*}
or
 \begin{align*}
 \brab{d^{(n)}_{0}-\varepsilon^{2}\Lambda_{h}}\bar{e}^n
=&\,\sum_{j=1}^{n-1}\brab{d^{(n)}_{n-j-1}-d^{(n)}_{n-j}-\varepsilon^{2}\eta^{n-j}\Lambda_{h}}\bar{e}^j\\
&\,+f(u^{n})-f(U^{n})+\Upsilon^{n}+R^{n},\quad{1}\leq{n}\leq{N}.
\end{align*}
By applying Lemma \ref{lem:Matrix-Inf-Norm} and the estimate \eqref{CONproof: nonlinear term}, one gets
\begin{align*}
d^{(n)}_{0}\mynormb{\bar{e}^n}_{\infty} &\le \,\mynormb{\brab{d^{(n)}_{0}-\varepsilon^{2}\Lambda_{h}}\bar{e}^n}_{\infty} \\
&\le\sum_{j=1}^{n-1}\mynormb{\brab{d^{(n)}_{n-j-1}-d^{(n)}_{n-j}-\varepsilon^{2}\eta^{n-j}\Lambda_{h}}\bar{e}^j}\\
&\quad +2\mynormb{e^{n}}_{\infty}+\mynormb{\Upsilon^{n}}_{\infty}+\mynormb{R^{n}}_{\infty},\quad{1}\leq{n}\leq{N}.
\end{align*}
Under the time-step constraintt \eqref{eq: tau condition-Max-Principle}, Lemma \ref{lem: Matrix-Qj-Inf-Norm}
with $S_n=0$ yields
\begin{align*}
\mynormb{\brab{d^{(n)}_{n-j-1}-d^{(n)}_{n-j}-\varepsilon^{2}\eta^{n-j}\Lambda_{h}}\bar{e}^j}_{\infty}
\le \brab{d^{(n)}_{n-j-1}-d^{(n)}_{n-j}}\mynormb{\bar{e}^j}_{\infty},\quad{1}\leq{j}\leq{n-1}.
\end{align*}
Thus, by applying the substitution formula \eqref{eq: substitution formulas}
and the triangle inequality,
it follows that
\begin{align*}
d^{(n)}_{0}\mynormb{\bar{e}^n}_{\infty}
\le\sum_{j=1}^{n-1}\brab{d^{(n)}_{n-j-1}-d^{(n)}_{n-j}}\mynormb{\bar{e}^j}_{\infty}
+2\sum_{j=1}^{n}\eta^{n-j}\mynormb{\bar{e}^j}_{\infty}
+\mynormb{\Upsilon^{n}}_{\infty}+\mynormb{R^{n}}_{\infty},
\end{align*}
or
\begin{align*}
\sum_{j=1}^{n}d^{(n)}_{n-j}\diff\mynormb{\bar{e}^j}_{\infty}\le2\sum_{j=1}^{n}\eta^{n-j}\mynormb{\bar{e}^j}_{\infty}
+\mynormb{\Upsilon^{n}}_{\infty}+\mynormb{R^{n}}_{\infty},\quad\text{$1\leq{n}\leq{N}$.}
\end{align*}
Under the choice \eqref{eq: ETA Choice}, one has $\eta\in [\frac12,1)$.
It is easy to check that the time-step {constraint} \eqref{eq: tau condition-Max-Principle} implies $\tau_n\le\frac{1+2r_n}{4(1+r_n)}$
or $b_0^{(n)}\ge4$. So Lemma \ref{lem: gronwall} with $\kappa=2$ and $\lambda:=\eta$ yields
\begin{align*}
\mynormb{\bar{e}^n}_{\infty}\le2\exp\left(\frac{4t_n}{1-\eta}\right)\sum_{j=1}^n\frac{1}{b^{(j)}_{0}}
\bra{\mynormb{\Upsilon^{j}}_{\infty}+\mynormb{R^{j}}_{\infty}}\quad\text{for $1\le n\le N$.}
\end{align*}
Then the substitution formula \eqref{eq: substitution formulas} gives
\begin{align}\label{CONproof: final estimate}
\mynormb{e^n}_{\infty}\le\frac{2}{1-\eta}\exp\left(\frac{4t_n}{1-\eta}\right)\sum_{j=1}^n\frac{1}{b^{(j)}_{0}}
\bra{\mynormb{\Upsilon^{j}}_{\infty}+\mynormb{R^{j}}_{\infty}}\quad\text{for $1\le n\le N$.}
\end{align}

Obviously, $\mynormb{R^{j}}_{\infty}\le C_uh^2$ for $j\ge1$ and thus we have
$$\sum_{j=1}^n\frac{1}{b^{(j)}_{0}}\mynormb{R^{j}}_{\infty}\le \sum_{j=1}^n\tau_j\mynormb{R^{j}}_{\infty}\le C_ut_nh^2.$$
By the Taylor's expansion (e.g., \cite[Theorem 10.5]{Thomee:2006}),
one has $\Upsilon^{1}=-\frac{1}{\tau_1}\int_{t_{0}}^{t_1}t\,\partial_{tt}u(t)\zd{t}$ and
\begin{align*}
\Upsilon^{n}=&\,-\frac{1+r_n}{2\tau_n}\int_{t_{n-1}}^{t_n}\!\!(t-t_{n-1})^2\partial_{ttt}u(t)\zd{t}
+\frac{r_n^2}{2(1+r_n)\tau_n}\int_{t_{n-2}}^{t_n}\!\!(t-t_{n-2})^2\partial_{ttt}u(t)\zd{t}, \,\, n \geq 2.
\end{align*}
We have
$\mynormb{\Upsilon^{1}}_{\infty}\le b^{(1)}_{0}\tau_1\int_{t_{0}}^{t_1}\mynormb{\partial_{tt}u(t)}_{\infty}\zd{t}$ and
\begin{align*}
\mynormb{\Upsilon^{j}}_{\infty}\le&\,\frac{1+r_j}{2}\tau_j\int_{t_{j-1}}^{t_j}\mynormb{\partial_{ttt}u(t)}_{\infty}\zd{t}
+\frac{r_j^2(\tau_j+\tau_{j-1})^2}{2(1+r_j)\tau_j}\int_{t_{j-2}}^{t_j}\mynormb{\partial_{ttt}u(t)}_{\infty}\zd{t}\\
=&\,(1+r_j)\tau_j\int_{t_{j-1}}^{t_j}\mynormb{\partial_{ttt}u(t)}_{\infty}\zd{t}
+\frac{\tau_j(1+r_j)}{2}\int_{t_{j-2}}^{t_{j-1}}\mynormb{\partial_{ttt}u(t)}_{\infty}\zd{t}\\
=&\,(b^{(j)}_{0}-b^{(j)}_{1})\tau_j^2\left(\int_{t_{j-1}}^{t_j}\mynormb{\partial_{ttt}u(t)}_{\infty}\zd{t}
+\frac{1}{2}\int_{t_{j-2}}^{t_{j-1}}\mynormb{\partial_{ttt}u(t)}_{\infty}\zd{t}\right)
\quad\text{for $j\ge2$},
\end{align*}
where $b^{(j)}_{0}-b^{(j)}_{1}=(1+r_j)/\tau_j$ from the definition \eqref{eq: BDF2 kernels} has been used.
It follows that
\begin{align*}
\sum_{j=1}^n\frac{1}{b^{(j)}_{0}}\mynormb{\Upsilon^{j}}_{\infty}
=&\,\tau_1\int_{t_{0}}^{t_1}\mynormb{\partial_{tt}u(t)}_{\infty}\zd{t}+\sum_{j=2}^n\frac{1}{d^{(j)}_{0}}\mynormb{\Upsilon^{j}}_{\infty}\\
\le&\,\tau_1\int_{t_{0}}^{t_1}\mynormb{\partial_{tt}u(t)}_{\infty}\zd{t}+
\sum_{j=2}^n\brab{1-b^{(j)}_{1}/b^{(j)}_{0}}\tau_j^2\int_{t_{j-1}}^{t_j}\mynormb{\partial_{ttt}u(t)}_{\infty}\zd{t}\\
&\,+\frac12\sum_{j=1}^{n-1}\brab{1-b^{(j+1)}_{1}/b^{(j+1)}_{0}}r_{j+1}^2\tau_{j}^2\int_{t_{j-1}}^{t_{j}}\mynormb{\partial_{ttt}u(t)}_{\infty}\zd{t}\\
\le&\,\tau_1\int_{t_{0}}^{t_1}\mynormb{\partial_{tt}u(t)}_{\infty}\zd{t}+
8\sum_{j=1}^n\tau_j^2\int_{t_{j-1}}^{t_j}\!\!\mynormb{\partial_{ttt}u(t)}_{\infty}\zd{t}
\quad\text{for $n\ge1$},
\end{align*}
where the step-ratio restriction \textbf{S0} was applied. Therefore we obtain from \eqref{CONproof: final estimate} that
\begin{align*}
\mynormb{e^n}_{\infty}\le\frac{2}{1-\eta}\exp{\left(\frac{4t_n}{1-\eta}\right)}
\left(\tau_1\int_{t_{0}}^{t_1}\!\!\mynormb{\partial_{tt}u(t)}_{\infty}\zd{t}
+8\sum_{j=1}^n\tau_j^2\int_{t_{j-1}}^{t_j}\!\!\mynormb{\partial_{ttt}u(t)}_{\infty}\zd{t}+C_ut_nh^2\right).
\end{align*}
This completes the proof.
\end{proof}

\section{Numerical implementations}
\setcounter{equation}{0}
In this section, we shall provide with some details on the numerical implementations and present several numerical examples.
For the nonlinear BDF2 scheme \eqref{eq: BDF2 vector scheme}, we shall perform a simple Newton-type iteration procedure at each time level with a tolerance $10^{-12}$. {Always we choose the solution at the previous level as the initial value of Newton iteration.} For more advanced nonlinear solvers, one can refer to \cite{WangWangWise:2010,ChenWangWangWise:2014,YanChenWangWise:2018}.

\subsection{Adaptive time-stepping strategy}

In simulating the phase field problems,  the temporal evolution of phase variables involve
multiple time scales,  such as the coarsening dynamics problems discussed
in Example \ref{Dynamics-Coarsening}, an initial random perturbation evolves on a fast
time scale while later dynamic coarsening evolves on a very slow time scale.
Therefore, adaptive time-stepping strategy is more practical to efficiently
resolve widely varying time scales and to significantly reduce the computational cost.
On the other hand,
one remarkable advantage of maximum norm stable scheme is that it can be easily
combined with an adaptive time strategy,
which adjusts the size of time step based on the accuracy requirement only.
In this paper, we use Algorithm \ref{Adaptive-Time-Strategy} which is motivated by \cite{Gomez2011Provably}
to choose adaptive time steps.

\begin{algorithm}
\caption{Adaptive time-stepping strategy}
\label{Adaptive-Time-Strategy}
\begin{algorithmic}[1]
\Require{Given $u^{n}$ and time step $\tau_{n}$}
\State Compute $u_{1}^{n+1}$ by using first-order  scheme with time step $\tau_{n}$.
\State Compute $u_{2}^{n+1}$ by using second-order  scheme with time step $\tau_{n}$.
\State Calculate $e_{n+1}=\|u_{2}^{n+1}-u_{1}^{n+1}\|/\|u_{2}^{n+1}\|$.
\If {$e_{n}<tol$}
\State Update time-step size $\tau_{n+1}\leftarrow\min\{\max\{\tau_{\min},\tau_{ada}\},\tau_{\max}\}$.
\Else
\State Recalculate with time-step size  $\tau_{n}\leftarrow\min\{\max\{\tau_{\min},\tau_{ada}\},\tau_{\max}\}$.
\State Goto 1
\EndIf
\end{algorithmic}
\end{algorithm}
The first-order  and second-order schemes used in Algorithm \ref{Adaptive-Time-Strategy} refer to
the backward Euler method and adaptive BDF2 scheme in this article, respectively.
 The adaptive time step $\tau_{ada}$ is given by
\begin{align*}
\tau_{ada}\bra{e,\tau}
=\rho\bra{\frac{tol}{e}}^{\frac{1}{2}}\tau_{cur},
\end{align*}
in which $\rho$ is a default safety coefficient, $tol$ is a reference tolerance,
$e$ is the relative error at each time level,
and $\tau_{cur}$ is the current time step.
In addition,
$\tau_{\max}$ and $\tau_{\min}$ are the predetermined maximum and minimum time steps.
In our computation, if not explicitly specified, we choose $\rho=0.6$, $tol=10^{-4}$,
$\tau_{\max}=0.1$ and {$\tau_{\min}=10^{-3}$}.

\subsection{Numerical examples}

\begin{example}\label{Accuracy-Test-BDF2}
To test the accuracy, we first consider
$\partial_{t}u
=\frac{1}{8\pi^2}\Delta u -f(u)+ g(\mathbf{x},t)$
for $\mathbf{x}\in(0,1)^{2}$ and $0<t<1$
such that it has an exact solution $u=\sin(2\pi x)\sin(2\pi y)\sin t$.
\end{example}

The numerical accuracy in time of BDF2 scheme is examined by
using the random mesh,
that is,  $\tau_{k}:=T\epsilon_{k}/S$ for $1\leq k\leq N$,
where $S=\sum_{k=1}^{N}\epsilon_{k}$ and $\epsilon_{k}\in(0,1)$ are random numbers.
The maximum norm error $e(N):=\max_{1\leq{n}\leq{N}}\|U^{n}-u^{n}\|_{\infty}$ is recorded in each run  and the experimental  order of convergence is computed by
$$\text{Order}\approx\frac{\log\bra{e(N)/e(2N)}}{\log\bra{\tau(N)/\tau(2N)}},$$
where $\tau(N)$ denotes the maximum time-step size for total $N$ subintervals.
We take the spatial grid points $M_1=1024$ in each direction such that the temporal error dominates the spatial error in each run
and solve the problem with $T=1$.
The numerical results are listed in Table \ref{BDF2-Time-Error},
where the number of step-ratio $r_k\ge1+\sqrt{2}$ is also listed in the fifth column.
It is somewhat surprising that the nonuniform BDF2 scheme
on random meshes maintains second-order accuracy even when there exists large step-ratios that do not satisfy the requirement $r_k < 1+\sqrt{2}$.

\begin{table}[htb!]
\begin{center}
\caption{Numerical accuracy of BDF2 scheme at time $T=1$.}\label{BDF2-Time-Error} \vspace*{0.3pt}
\def\temptablewidth{0.6\textwidth}
{\rule{\temptablewidth}{0.5pt}}
\begin{tabular*}{\temptablewidth}{@{\extracolsep{\fill}}ccccc}
  $N$      &$\tau$       &$e(N)$     &Order   &$r_k\ge1+\sqrt{2}$ \\
\midrule
  10       &1.88e-01	 &2.56e-03   &$-$     &1\\
  20       &1.10e-01	 &8.16e-04	 &2.12    &4\\
  40       &4.67e-02	 &1.39e-04	 &2.06    &3\\
  80       &2.42e-02	 &3.41e-05	 &2.14    &9\\
\end{tabular*}
{\rule{\temptablewidth}{0.5pt}}
\end{center}
\end{table}	

\begin{figure}[htb!]
\centering
\includegraphics[width=3.0in,height=2.0in]{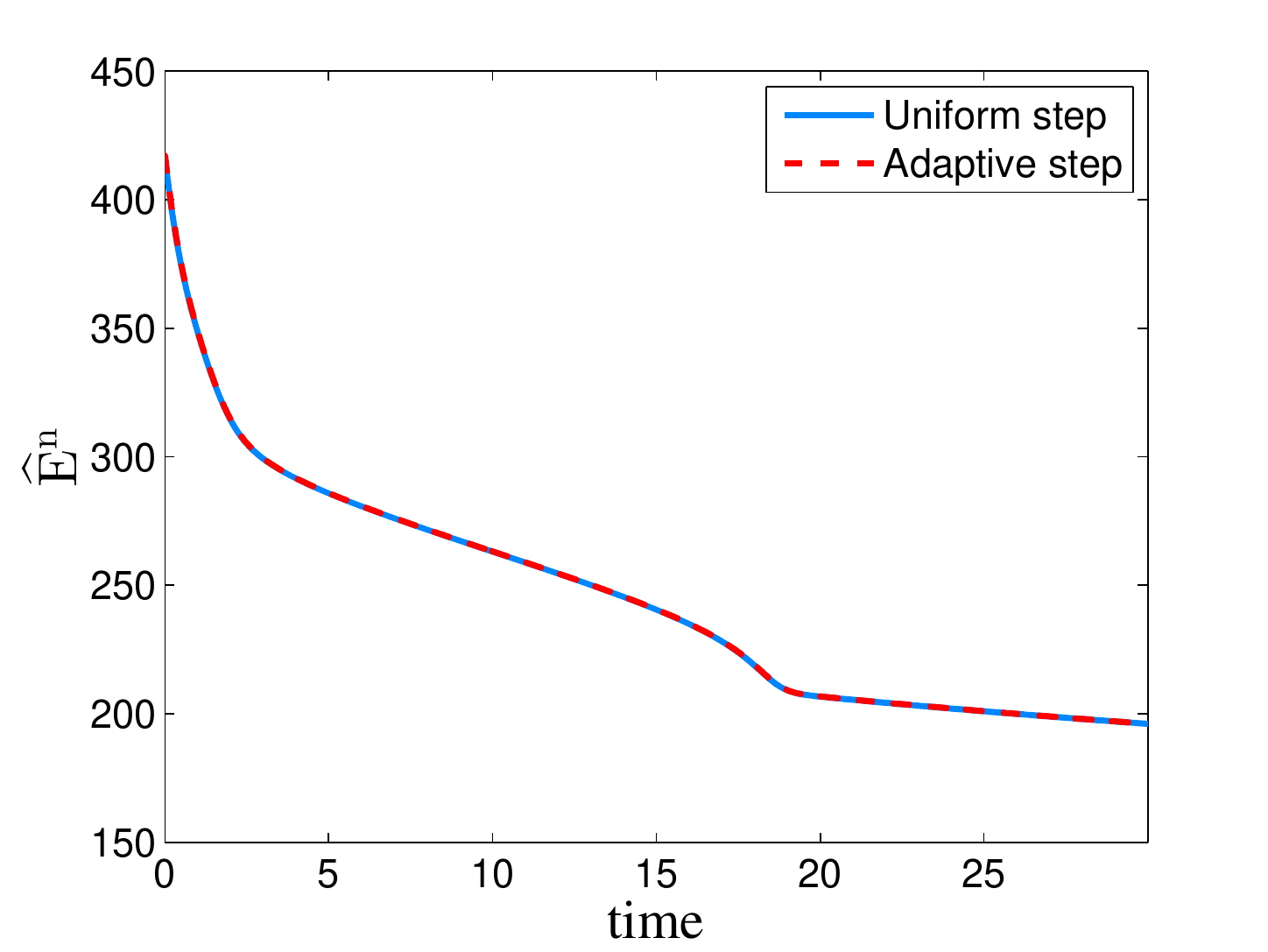}
\includegraphics[width=3.0in,height=2.0in]{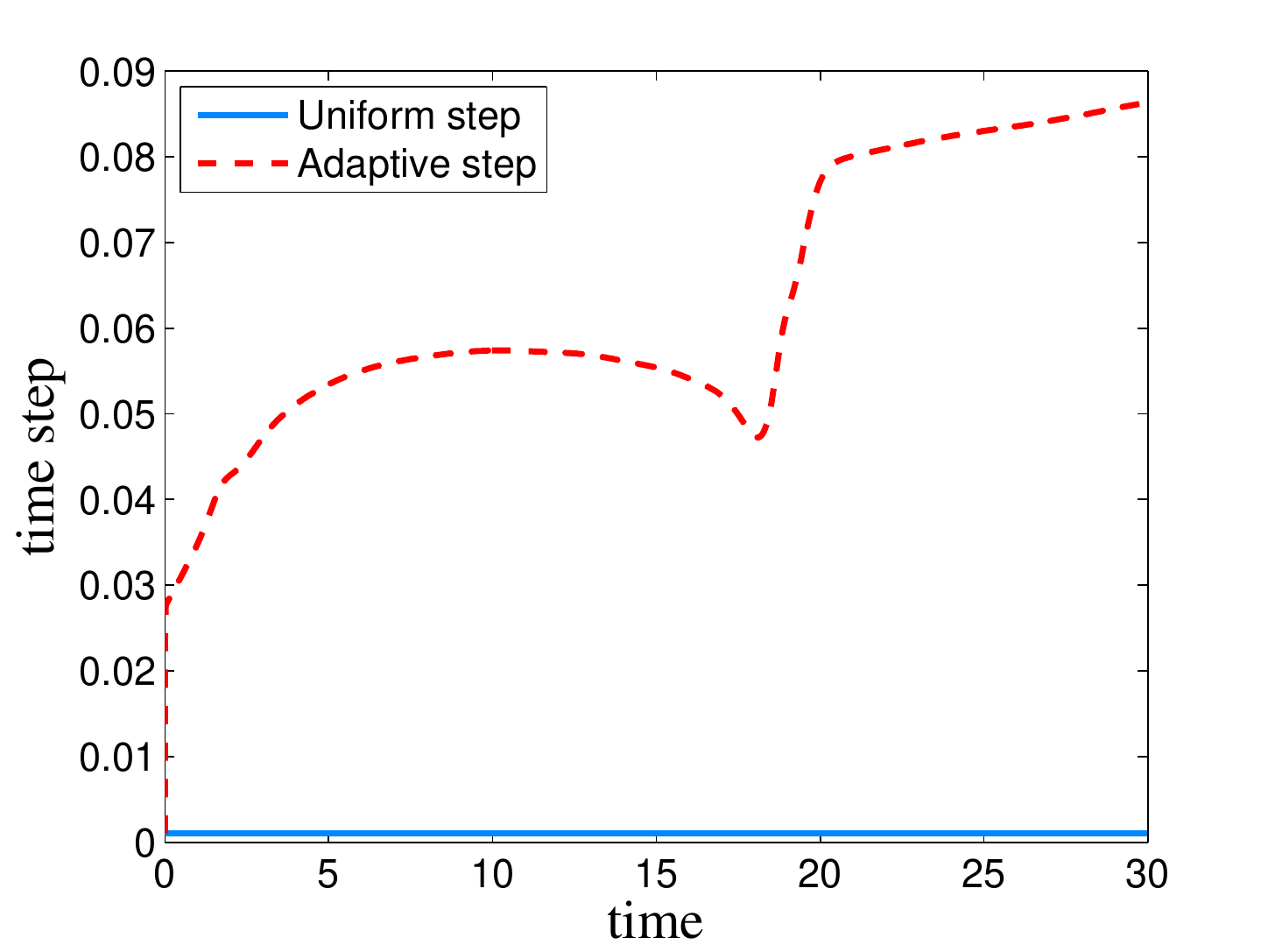}
\caption{Evolutions of energy (left) and time steps (right) of
   the  Allen-Cahn equation using different time strategies until time $T=30$.}
\label{Comparison-Uniform-Adaptive-Energy}
\end{figure}

\begin{example}\label{Simulating-Bubbles}
We next consider the Allen-Cahn model (\ref{Problem-1})-(\ref{Problem-2}) with the diffusion coefficient $\varepsilon=0.02$.
The nonuniform BDF2 scheme
is applied to simulate the merging of four bubbles
with an initial condition
\begin{align}
\phi_{0}\bra{\mathbf{x}}
=&-\tanh\bra{\bra{(x-0.3)^{2}+y^{2}-0.2^2}/\varepsilon}
\tanh\bra{\bra{(x+0.3)^{2}+y^{2}-0.2^2}/\varepsilon}\nonumber\\
&\times\tanh\bra{\bra{x^{2}+(y-0.3)^{2}-0.2^2}/\varepsilon}
\tanh\bra{\bra{x^{2}+(y+0.3)^{2}-0.2^2}/\varepsilon}.
\end{align}
The computational domain $\Omega=(-1,1)^{2}$ is divided uniformly into 128 parts in each direction.
\end{example}
We now examine different time strategies, i.e., the uniform and adaptive time approaches, for simulating the merging of four bubbles.
We start with the calculation of the solution  until the time $T=30$ with a constant time step $\tau=10^{-3}$.
We then implement the adaptive strategy described in Algorithm \ref{Adaptive-Time-Strategy}
to simulate the merging of bubbles.
The time evolution of discrete energies  and time steps  are depicted in Figure \ref{Comparison-Uniform-Adaptive-Energy}.
As can be seen, the  adaptive energy curve is practically indistinguishable  from
the one obtained using the small constant time step $\tau=10^{-3}$.
As a consequence, the total number of adaptive time steps are 511 while it takes 30000 steps for uniform
grid, showing that the time-stepping adaptive strategy is computationally efficient.

\begin{figure}[htb!]
\centering
\includegraphics[width=2.0in]{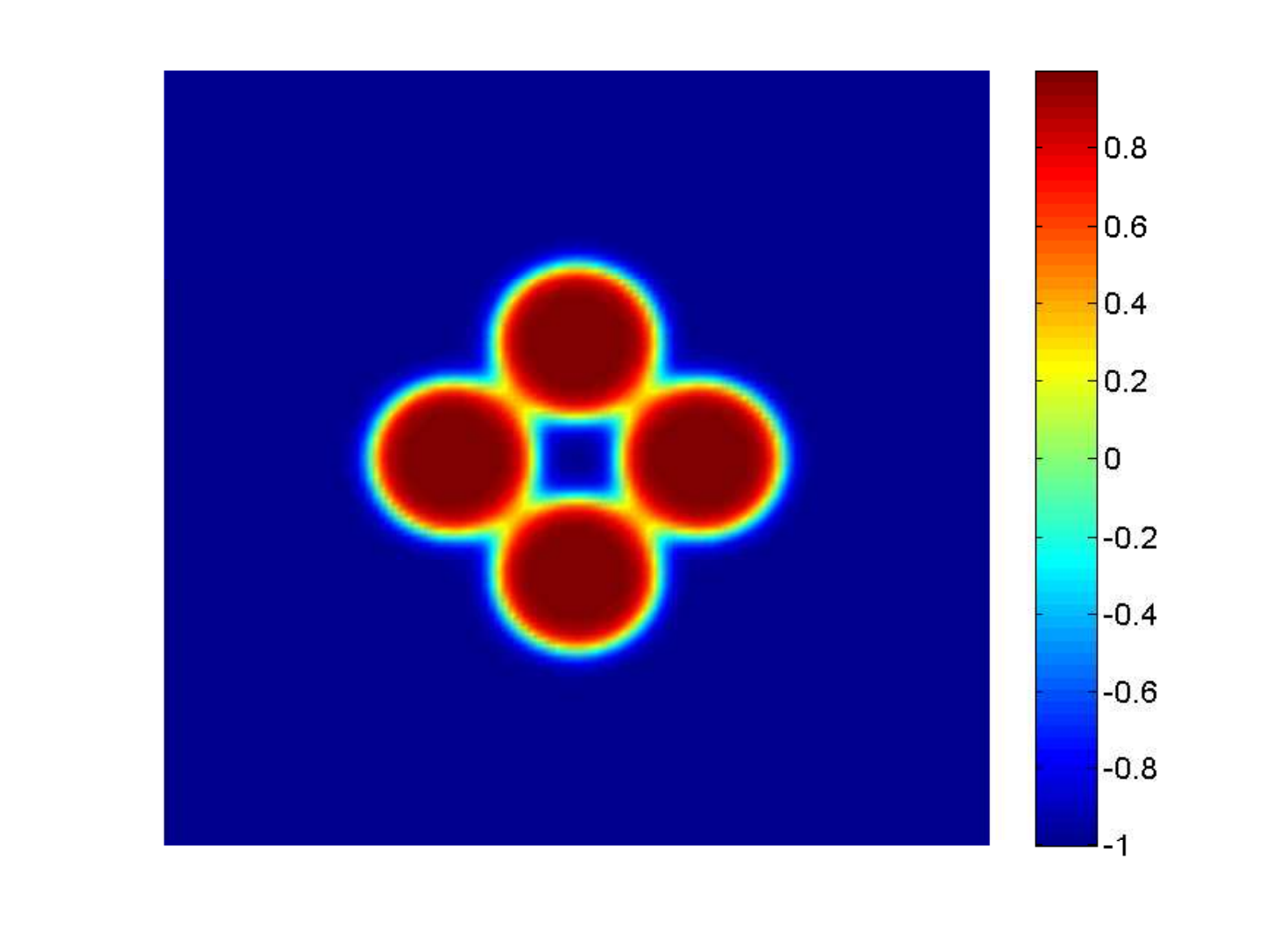}
\includegraphics[width=2.0in]{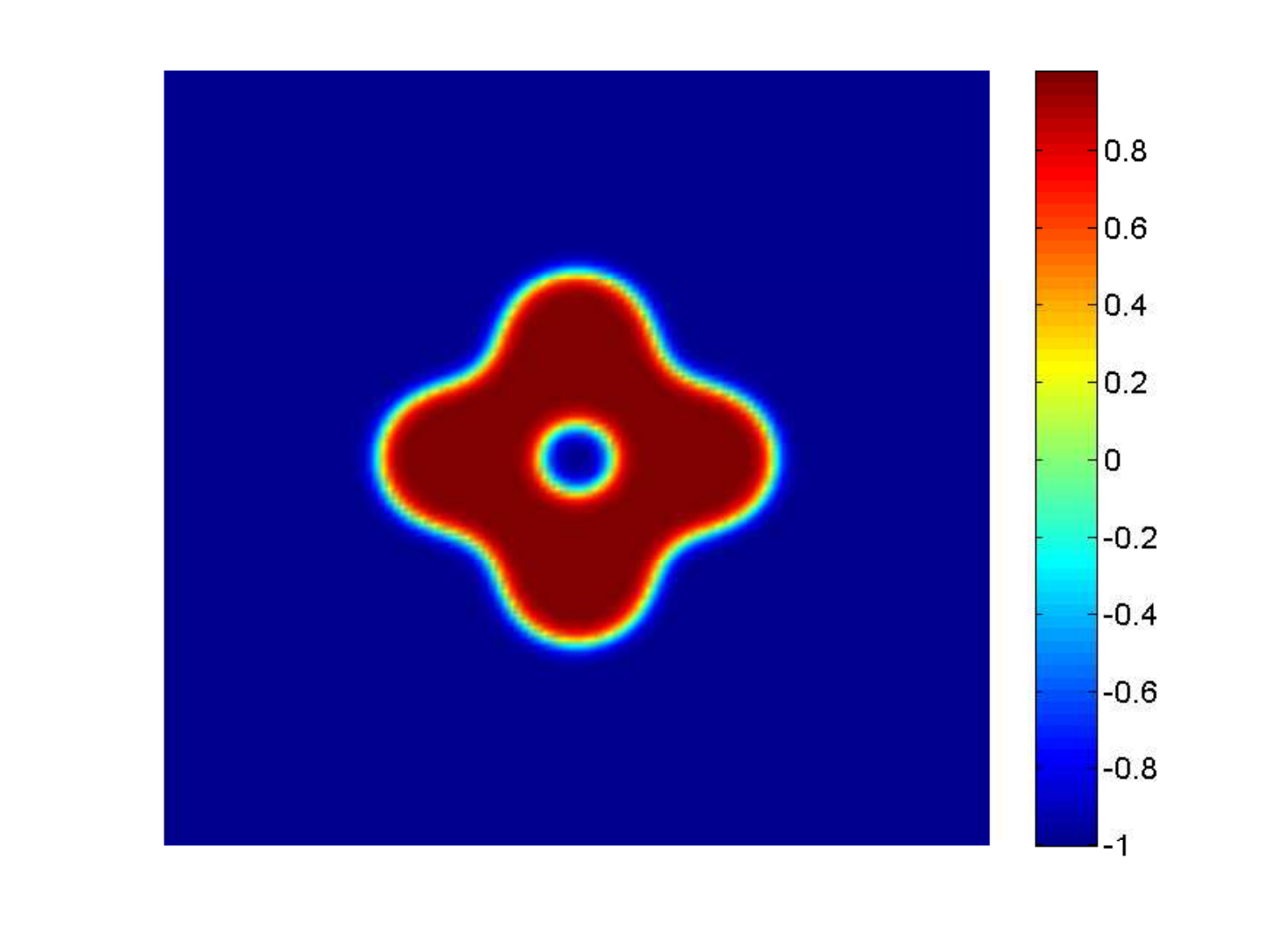}
\includegraphics[width=2.0in]{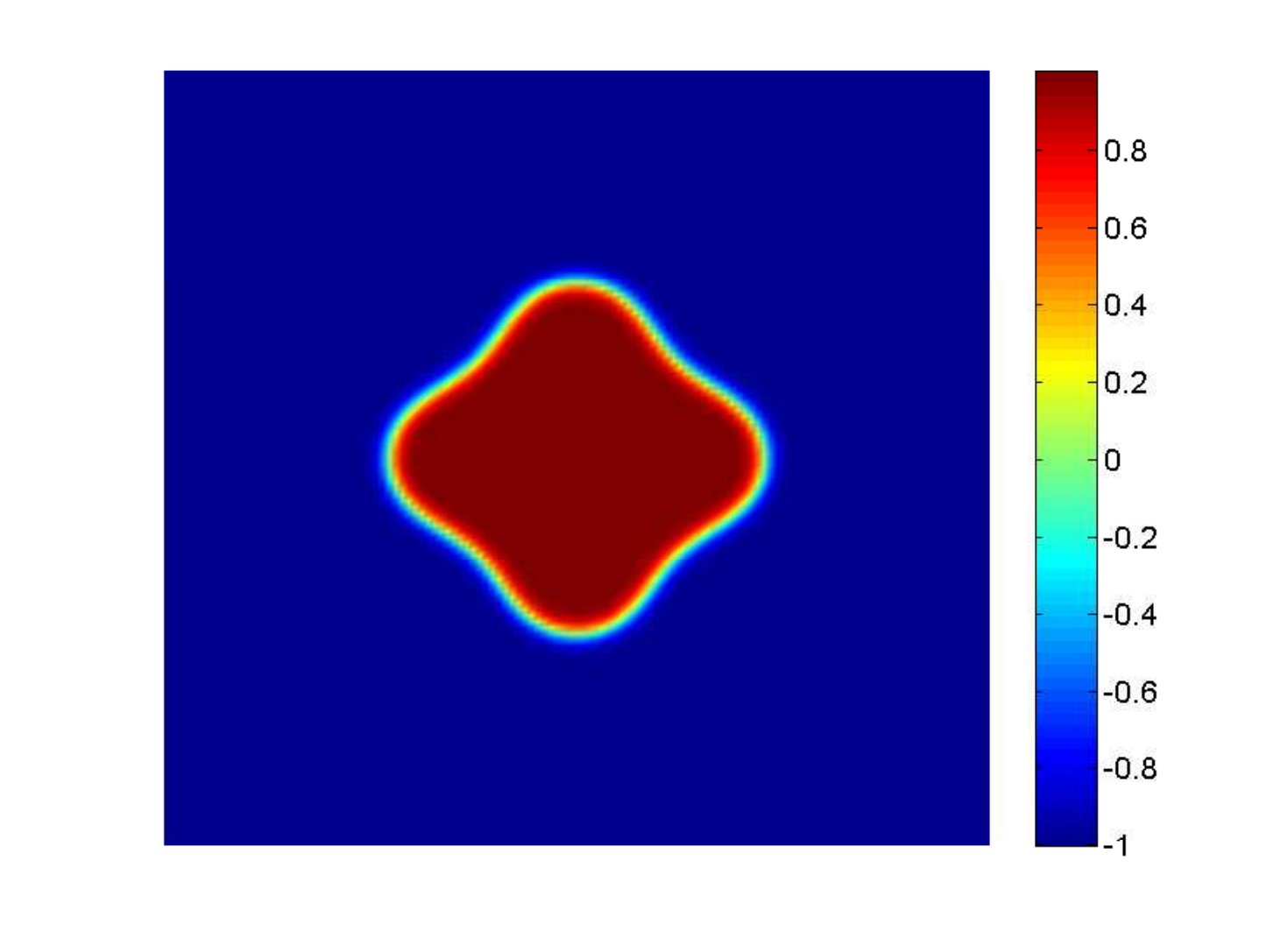}\\
\includegraphics[width=2.0in]{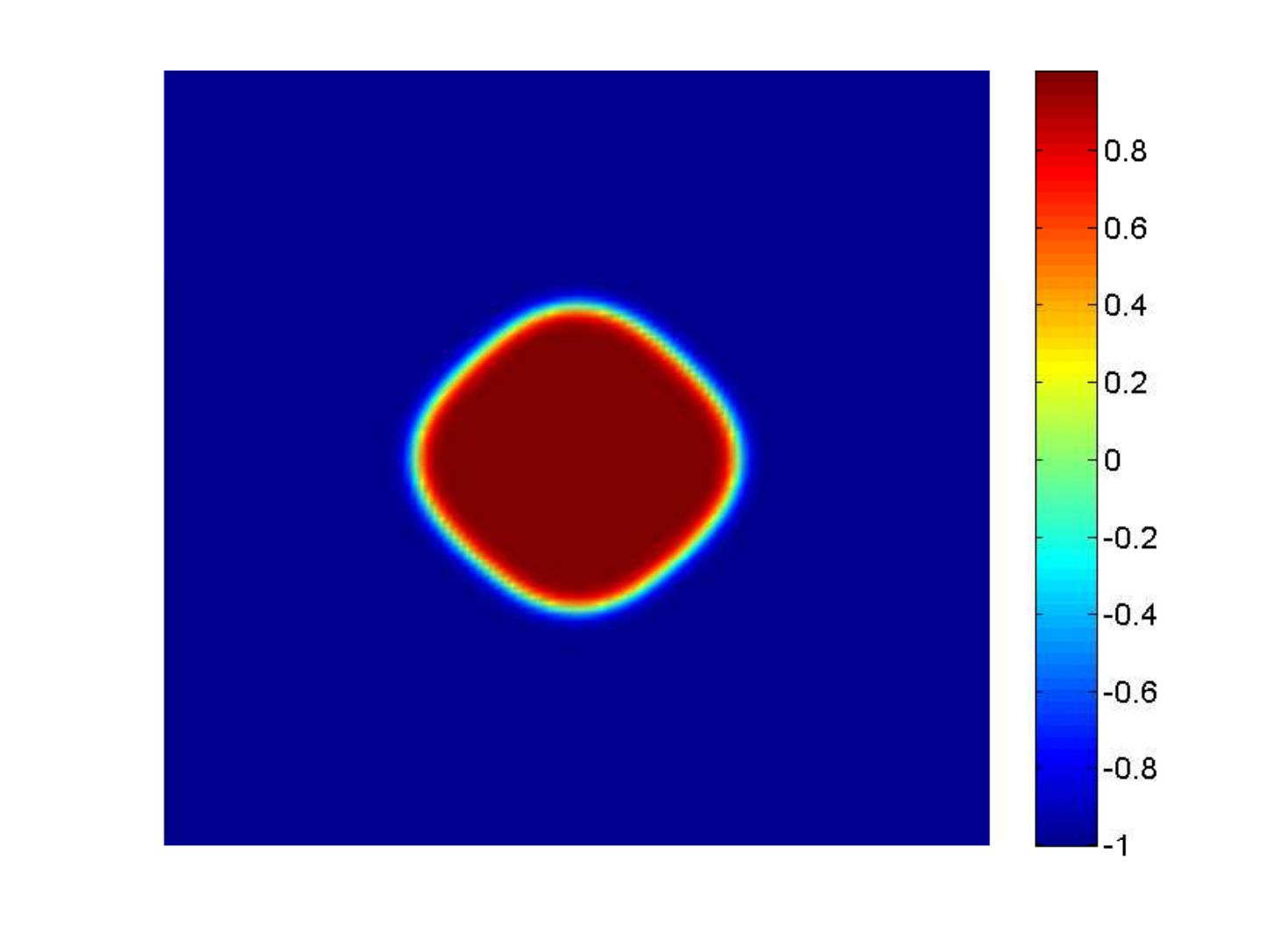}
\includegraphics[width=2.0in]{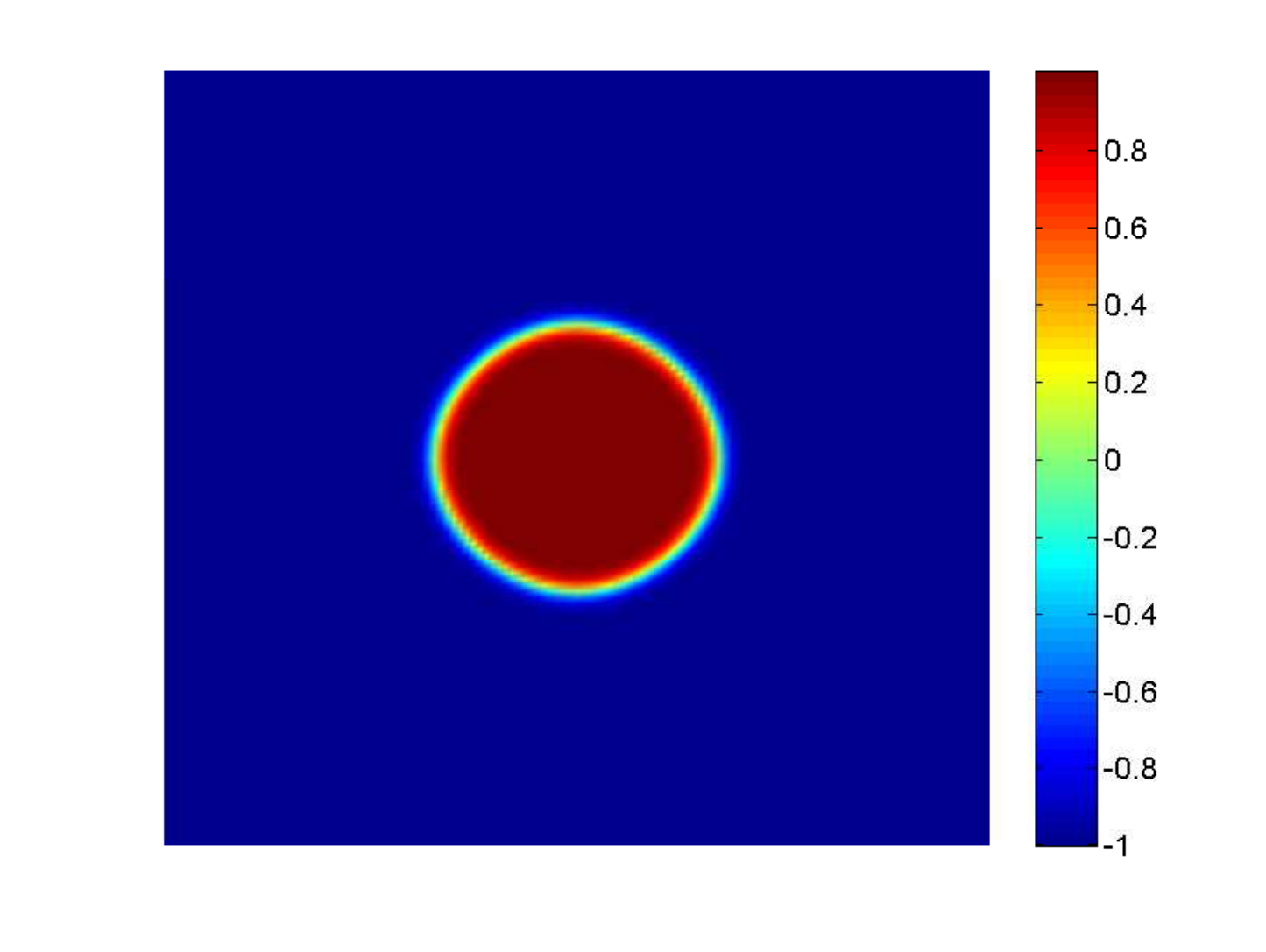}
\includegraphics[width=2.0in]{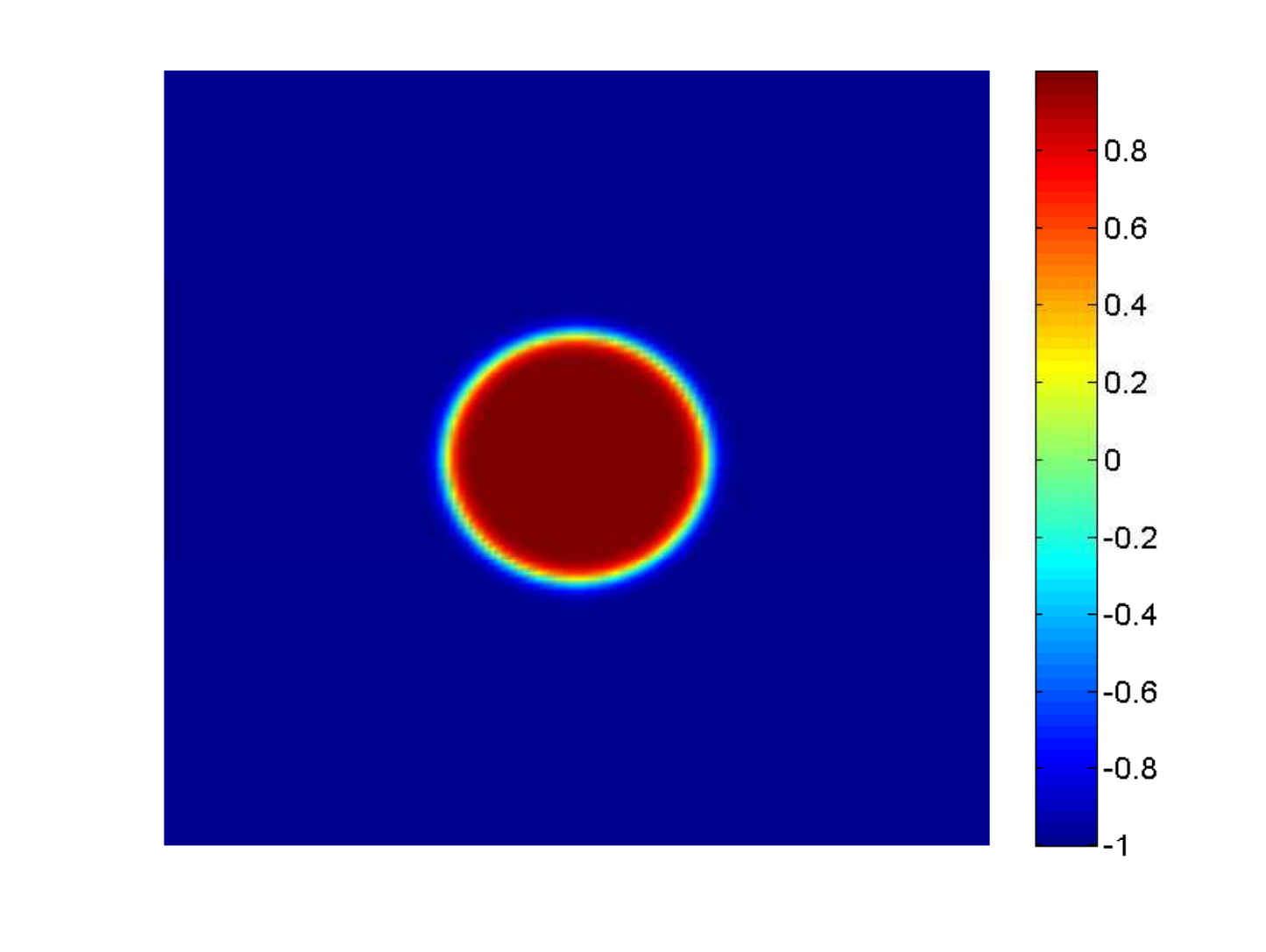}\\
\caption{Solution snapshots of the  Allen-Cahn equation using adaptive time strategy at
  $t=1, 10, 20, 50,80,100$, respectively.}
\label{Snapshots-Bubbles}
\end{figure}

\begin{figure}[htb!]
\centering
\includegraphics[width=3.0in,height=2.0in]{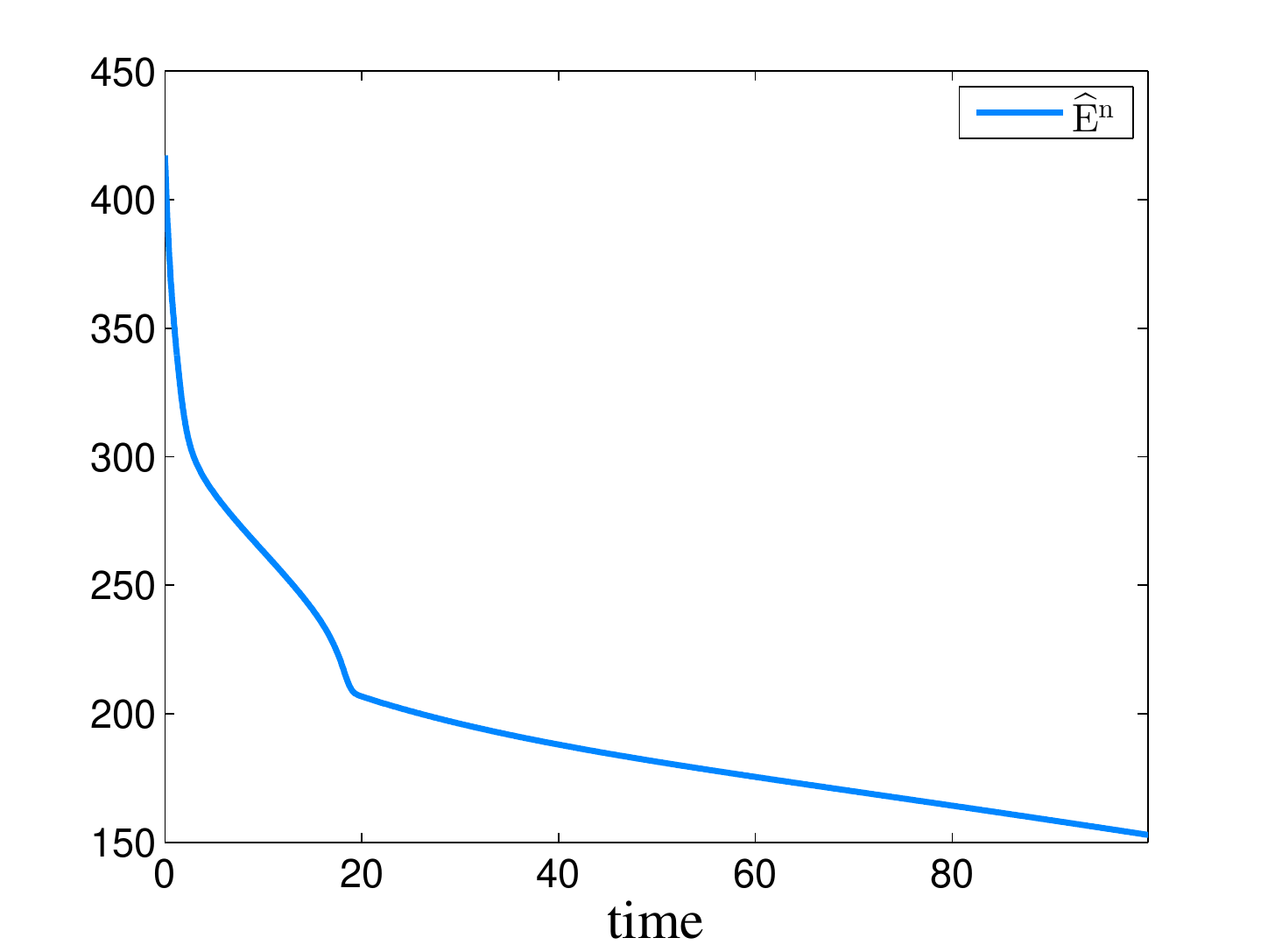}
\includegraphics[width=3.0in,height=2.0in]{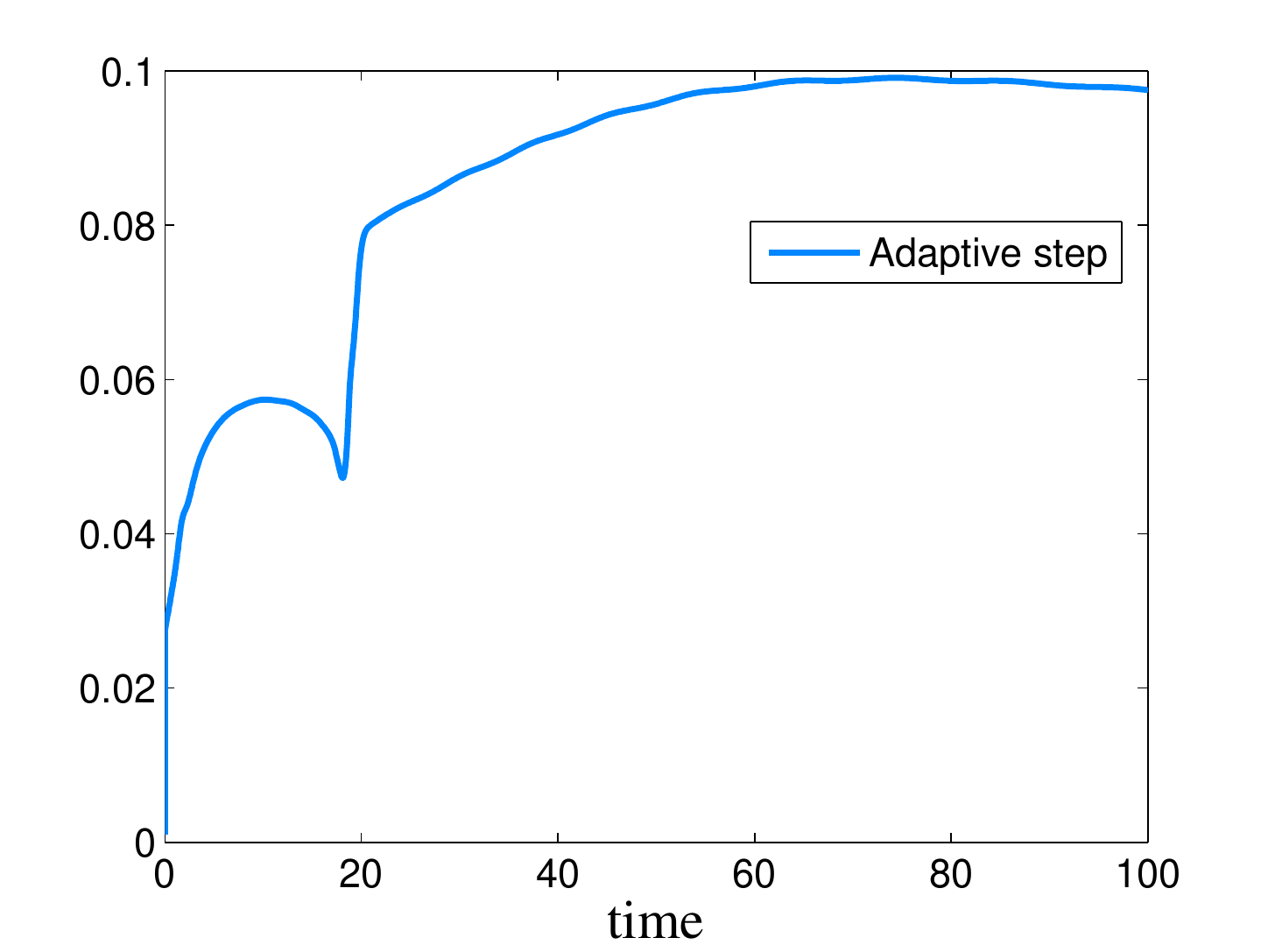}
\caption{Evolutions of energy (left) and time steps (right) of
   the  Allen-Cahn equation using adaptive time strategy until final time $T=100$.}
\label{Bubbles-Energy}
\end{figure}

We now apply the nonuniform BDF2 scheme coupled with the adaptive time strategy to
simulate the merging of bubbles with $T=100$.
The time evolution of the phase variable is summarized in Figure \ref{Snapshots-Bubbles}.
As can be seen in the figures,
the initial separated four bubbles gradually coalesce into a single big
bubble while the volume becomes smaller with time owing to that the Allen-Cahn model dose
not conserve the initial volume.
The discrete energy and adaptive time step are shown in Figure \ref{Bubbles-Energy}.
We observe that the energy evolution undergoes large variations initially and at time $t=20$,
but changes very little in other time intervals.
As a result, we see that small time steps are used when the energy variation is large,
while large time steps are utilized when the energy variation is small.

\begin{example}\label{Dynamics-Coarsening}
We finally consider the coarsening dynamics of the Allen-Cahn model
with the model parameter $\varepsilon=0.01$.
We choose a random initial condition $u_{0}=0.95+ {rand}(\mathbf{x})\times 0.05$ by assigning a random number varying from
$-0.05$ to $0.05$ at each grid points.
In the following computation, we use $128 \times 128$ uniform meshes in space
to discretize the domain $\Omega=(0,1)^{2}$.
\end{example}
\begin{figure}[htb!]
\centering
\includegraphics[width=2.0in]{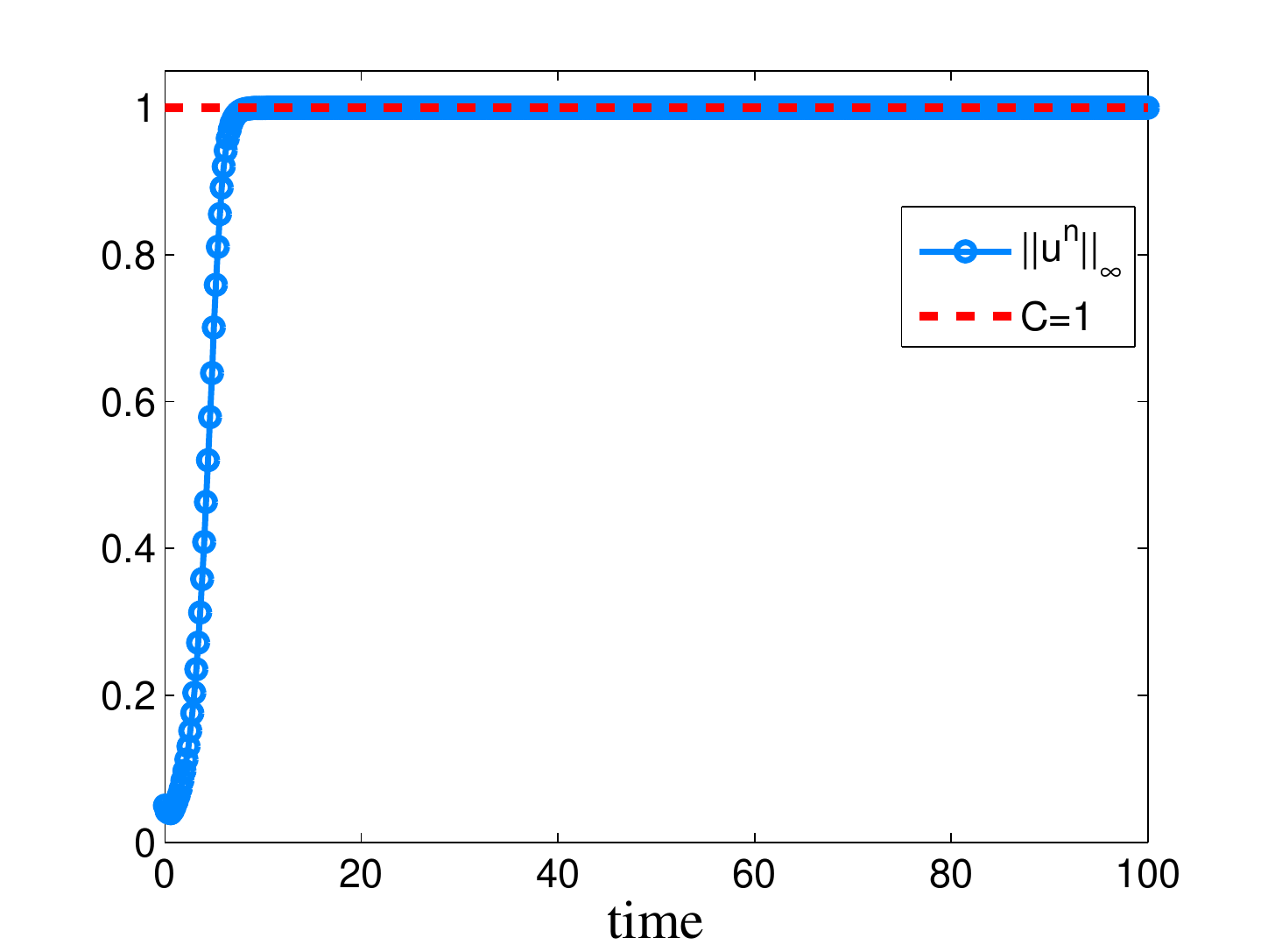}
\includegraphics[width=2.0in]{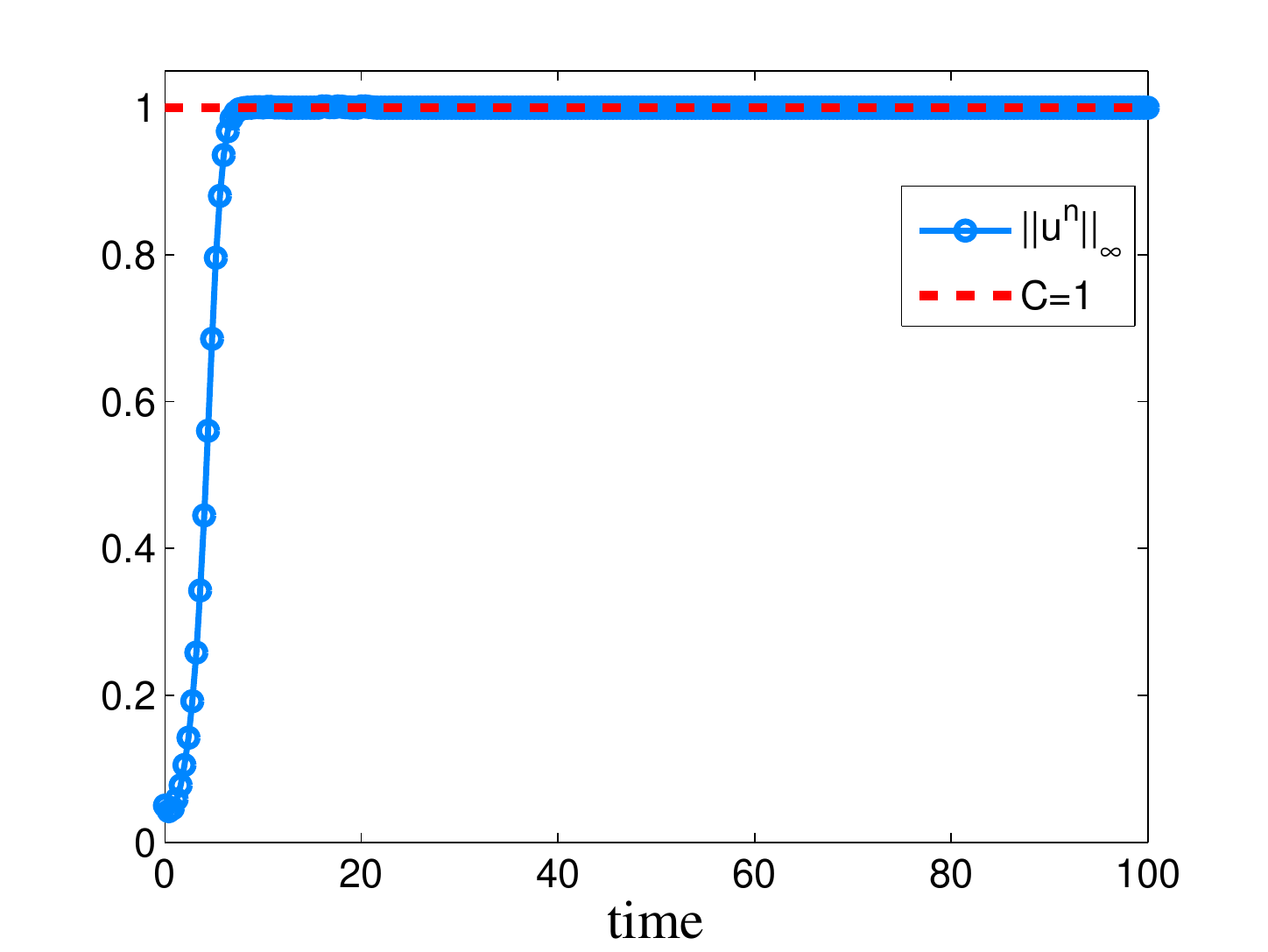}
\includegraphics[width=2.0in]{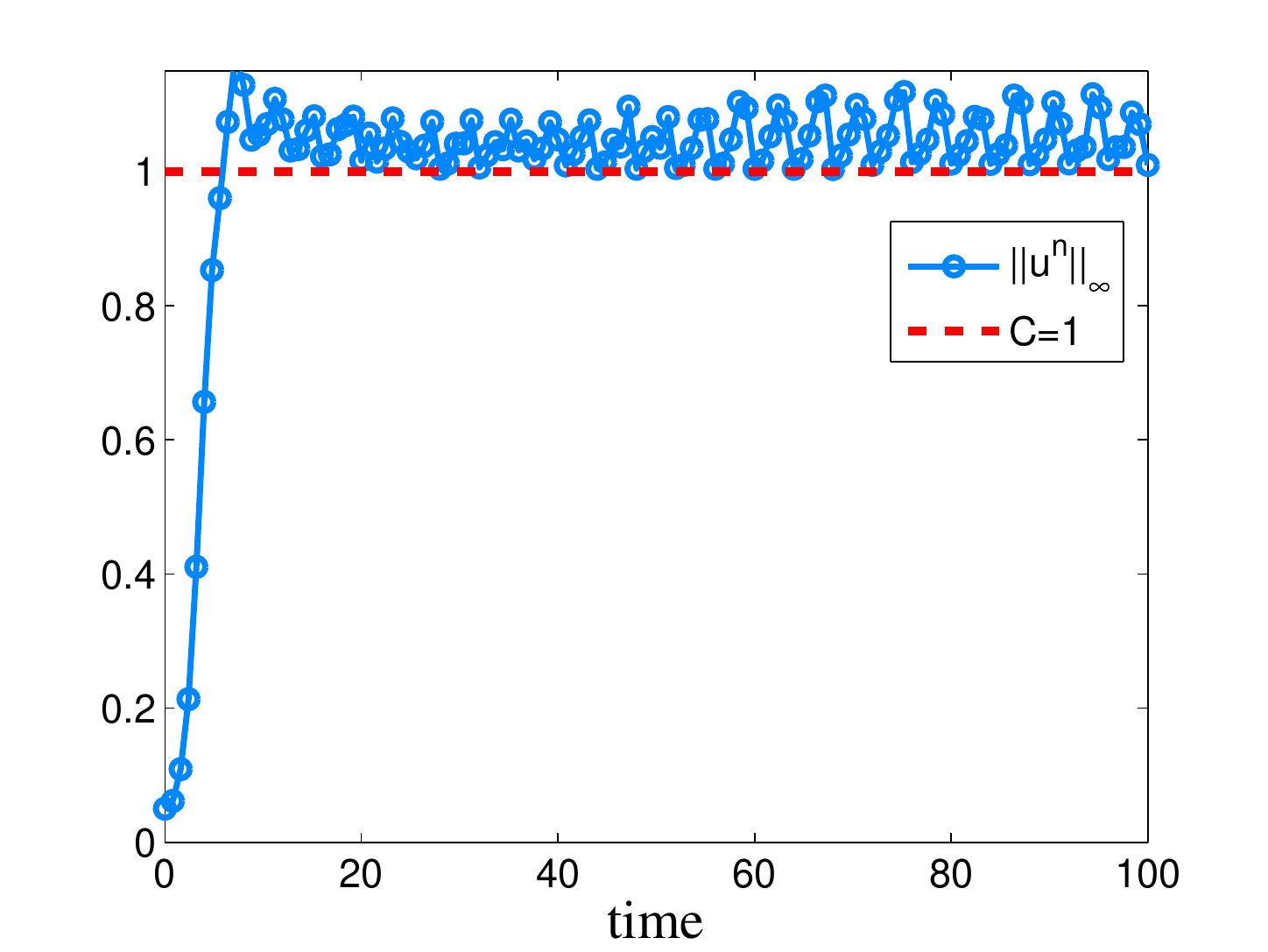}\\
\includegraphics[width=2.0in]{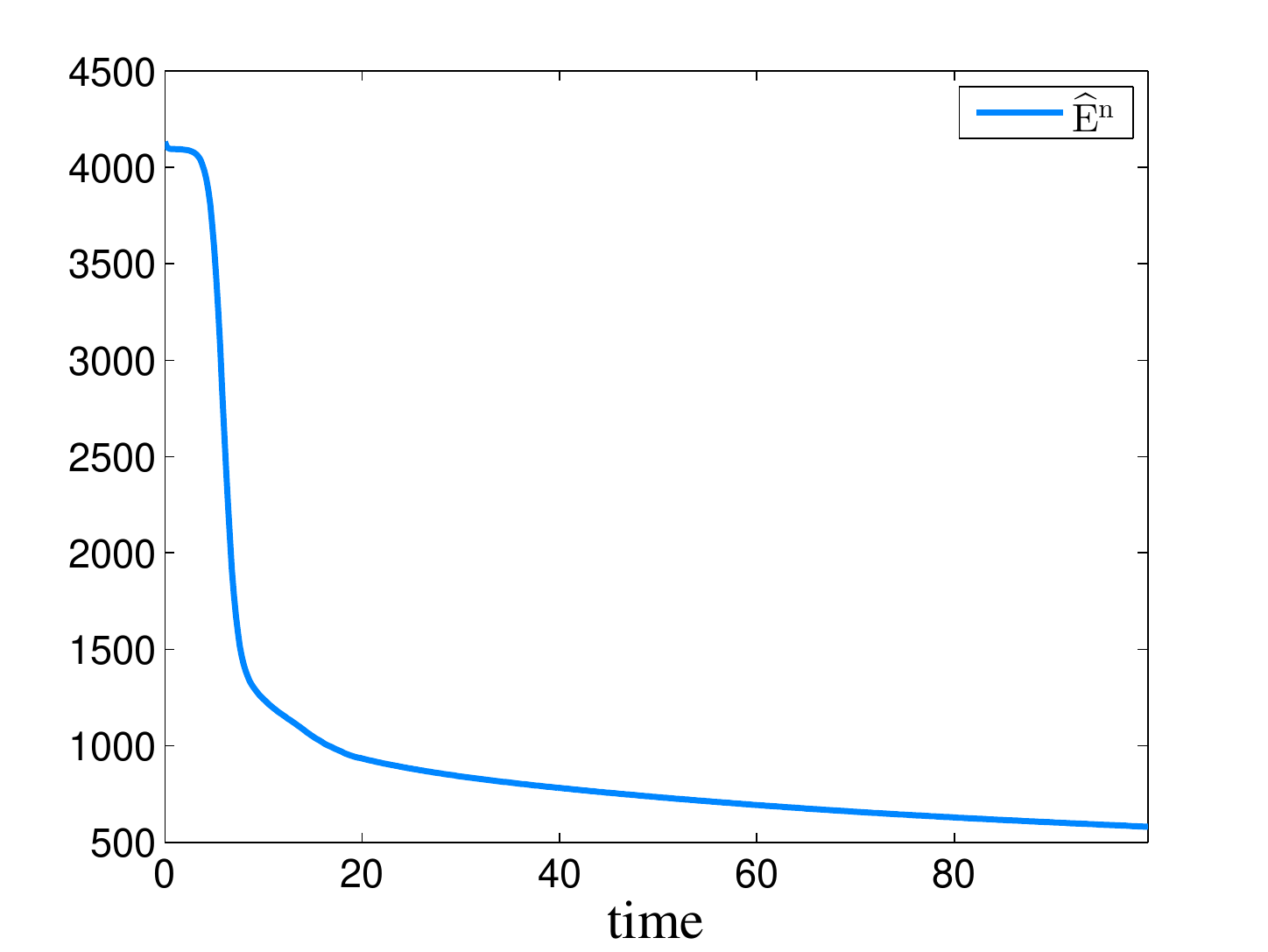}
\includegraphics[width=2.0in]{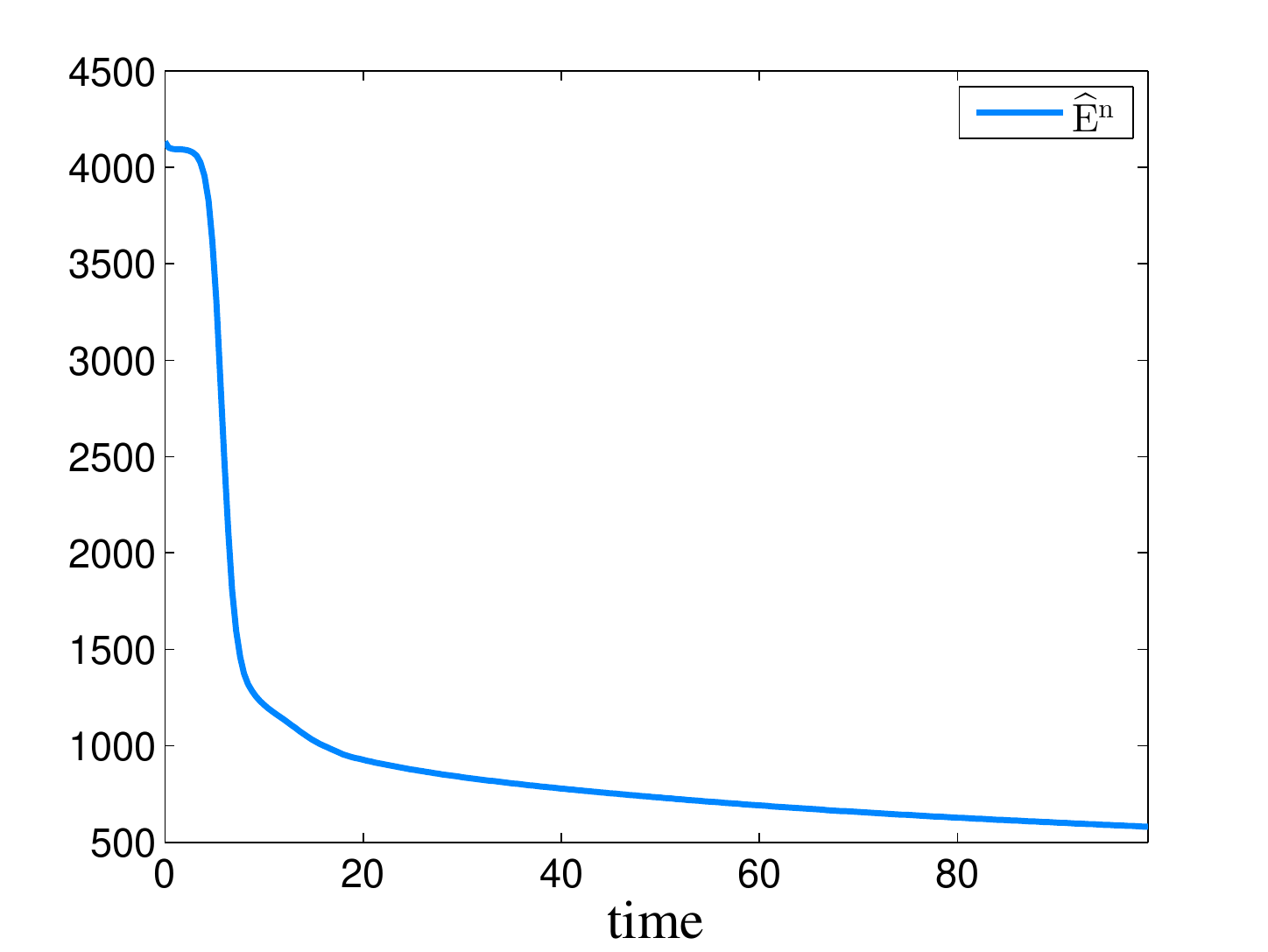}
\includegraphics[width=2.0in]{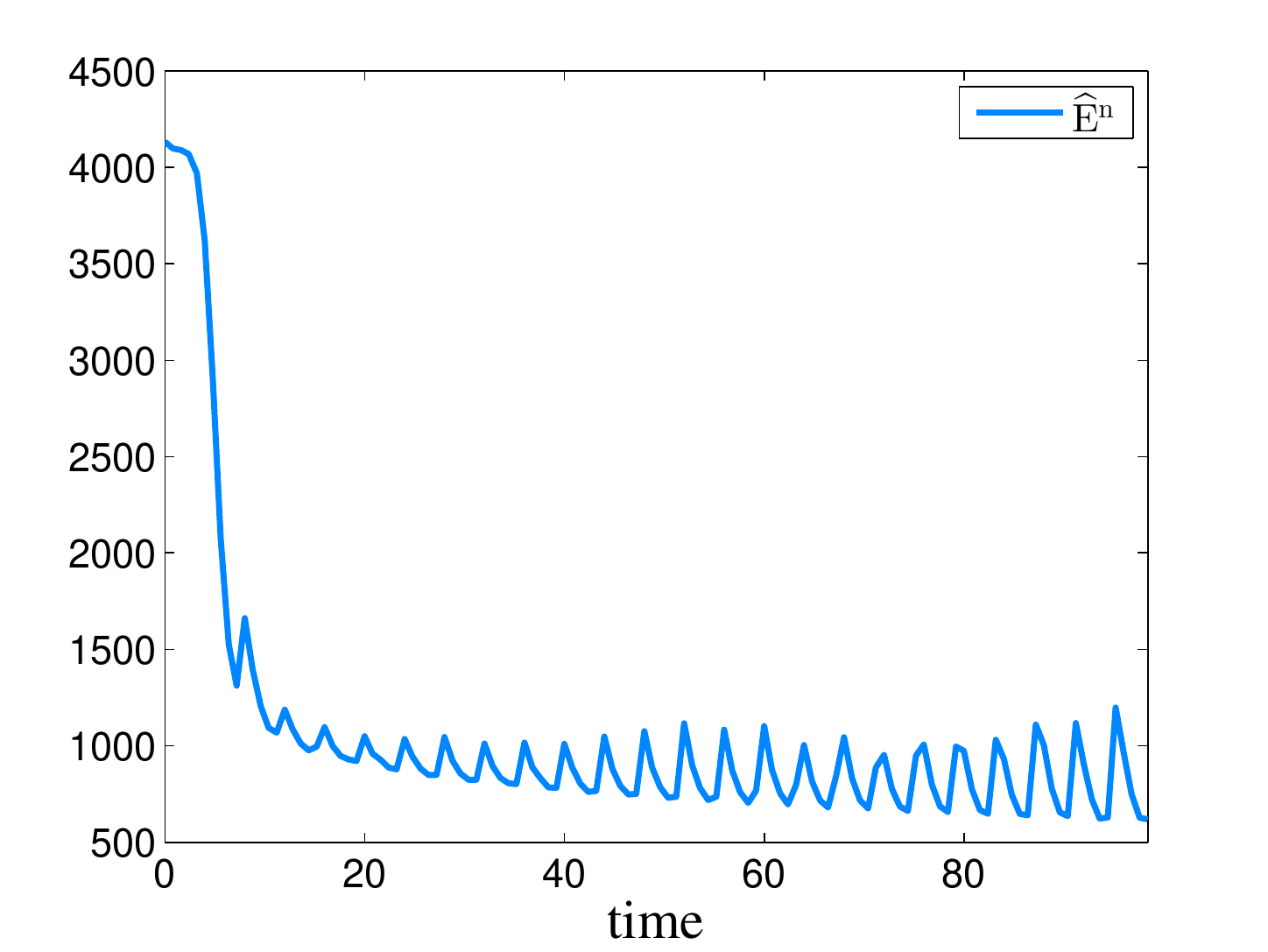}\\
\caption{Maximum norm (top) and energy (bottom) of the  Allen-Cahn equation using different time
  steps $\tau=0.2,0.4,0.8$ (from left to right), respectively.}
\label{MaxPrinciple-Different-Steps}
\end{figure}
We first investigate the effect of uniform time step size on the maximum norm and discrete energy.
The numerical results obtained from different time steps $\tau=0.2,0.4,0.8$ with $T=100$
are shown in Figure  \ref{MaxPrinciple-Different-Steps}.
As can be seen from the figures,
the maximum values of the numerical solutions are bounded by $1$ and
the energy dissipation law holds if time steps $\tau=0.2,0.4.$
These numerical results imply that the  constraintt \eqref{eq: tau condition-Max-Principle} for time step size
to ensure the discrete maximum principle is a sufficient condition.

\begin{figure}[htb!]
\centering
\includegraphics[width=2.0in]{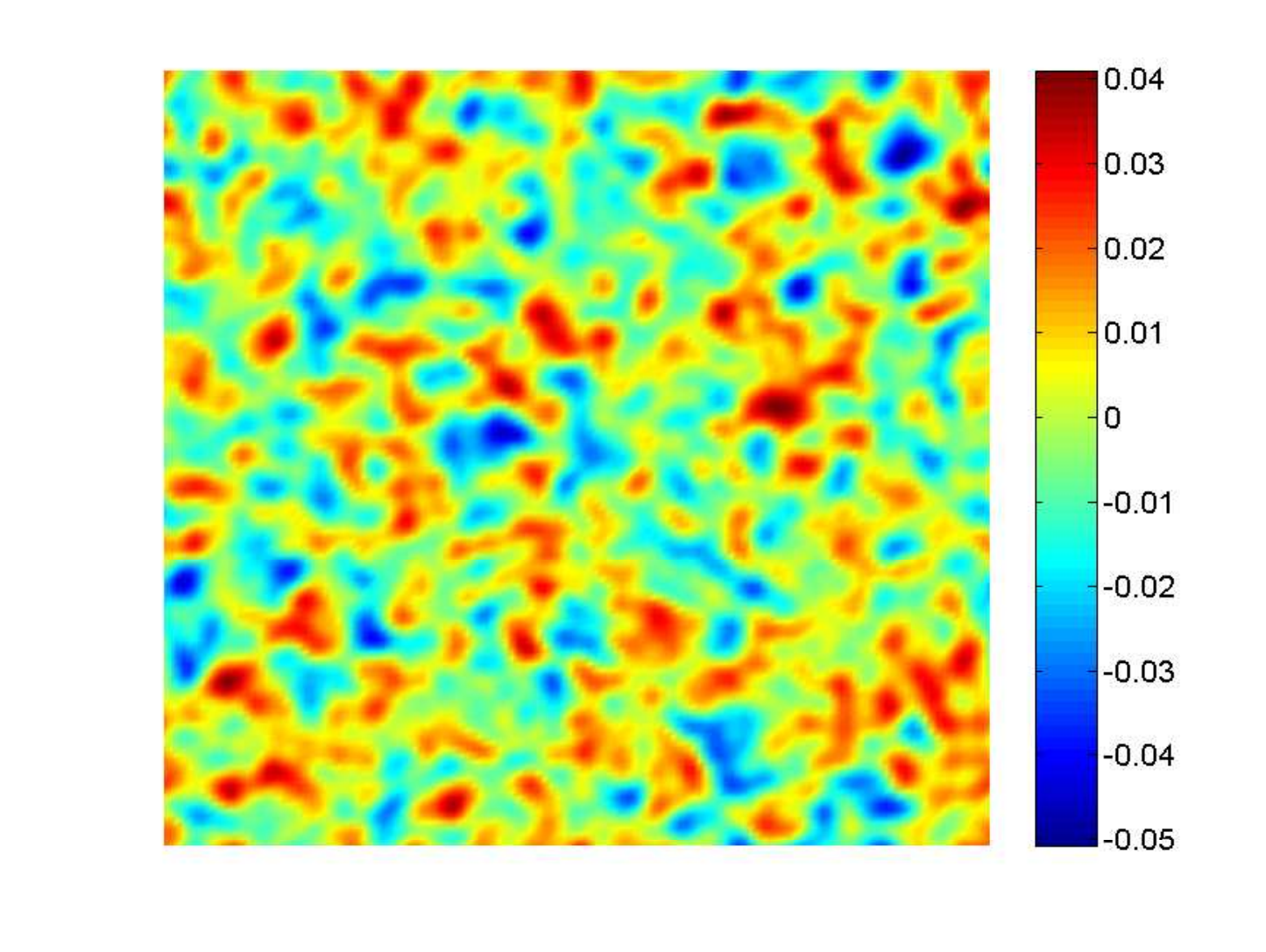}
\includegraphics[width=2.0in]{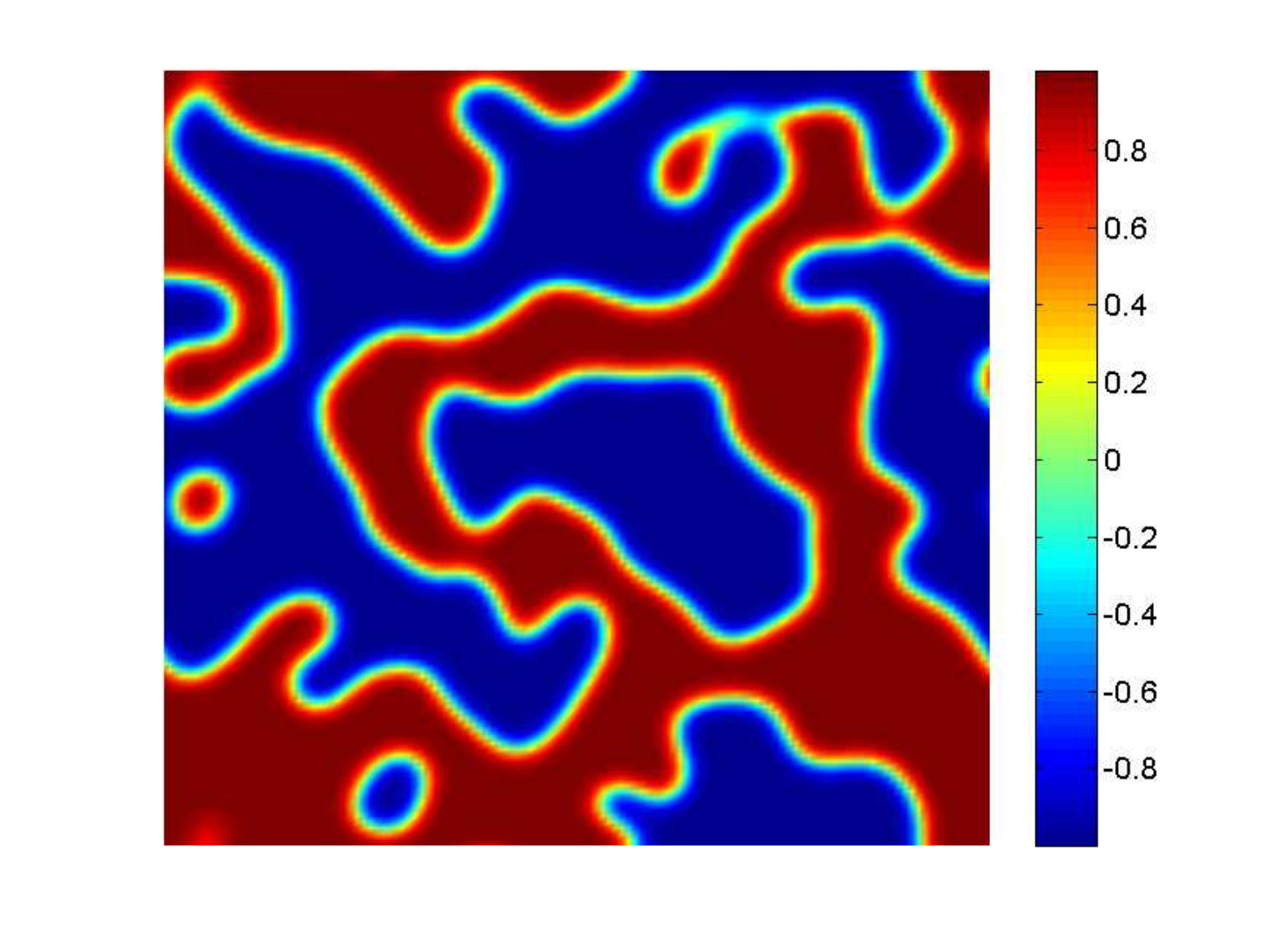}
\includegraphics[width=2.0in]{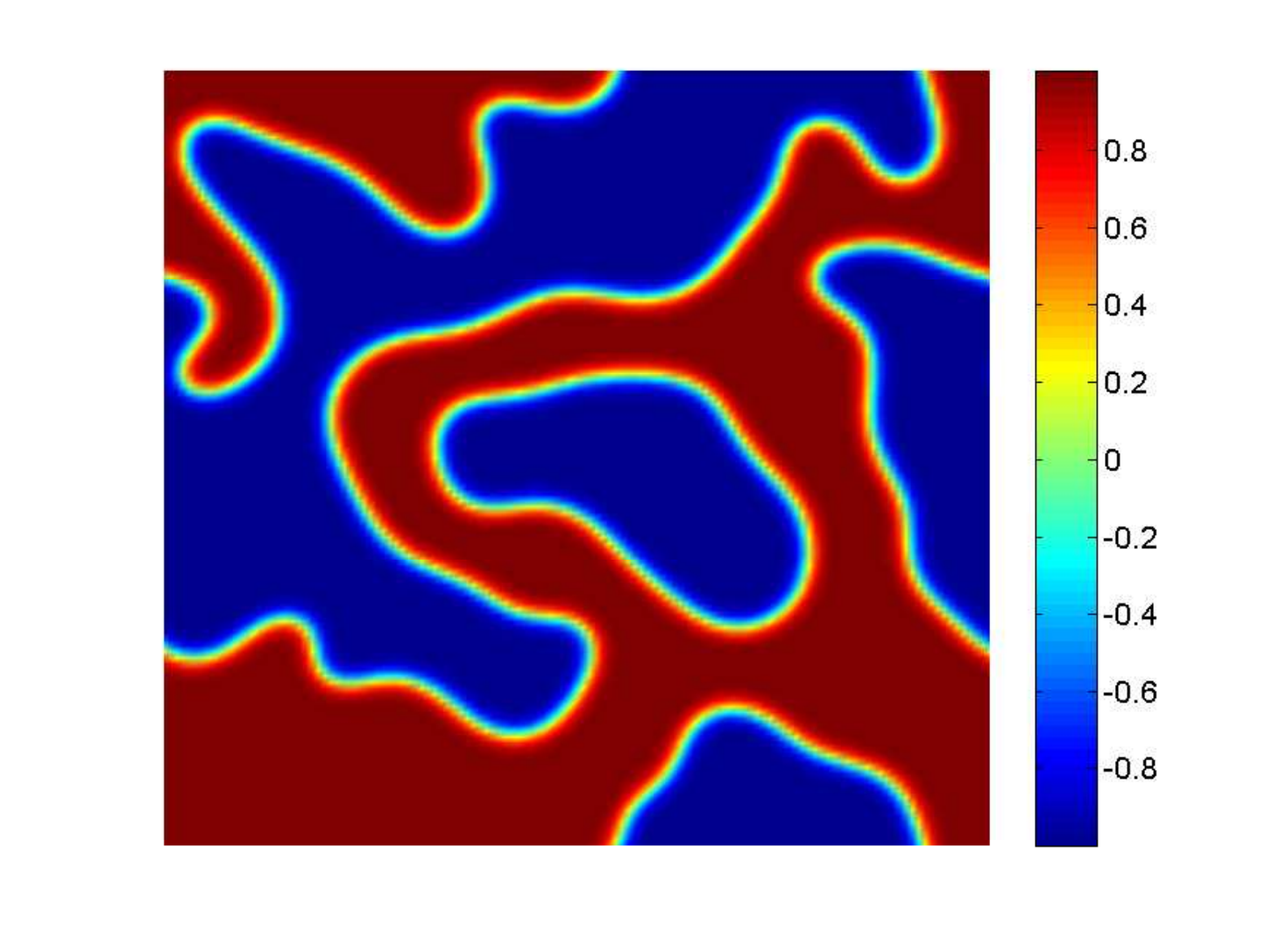}\\
\includegraphics[width=2.0in]{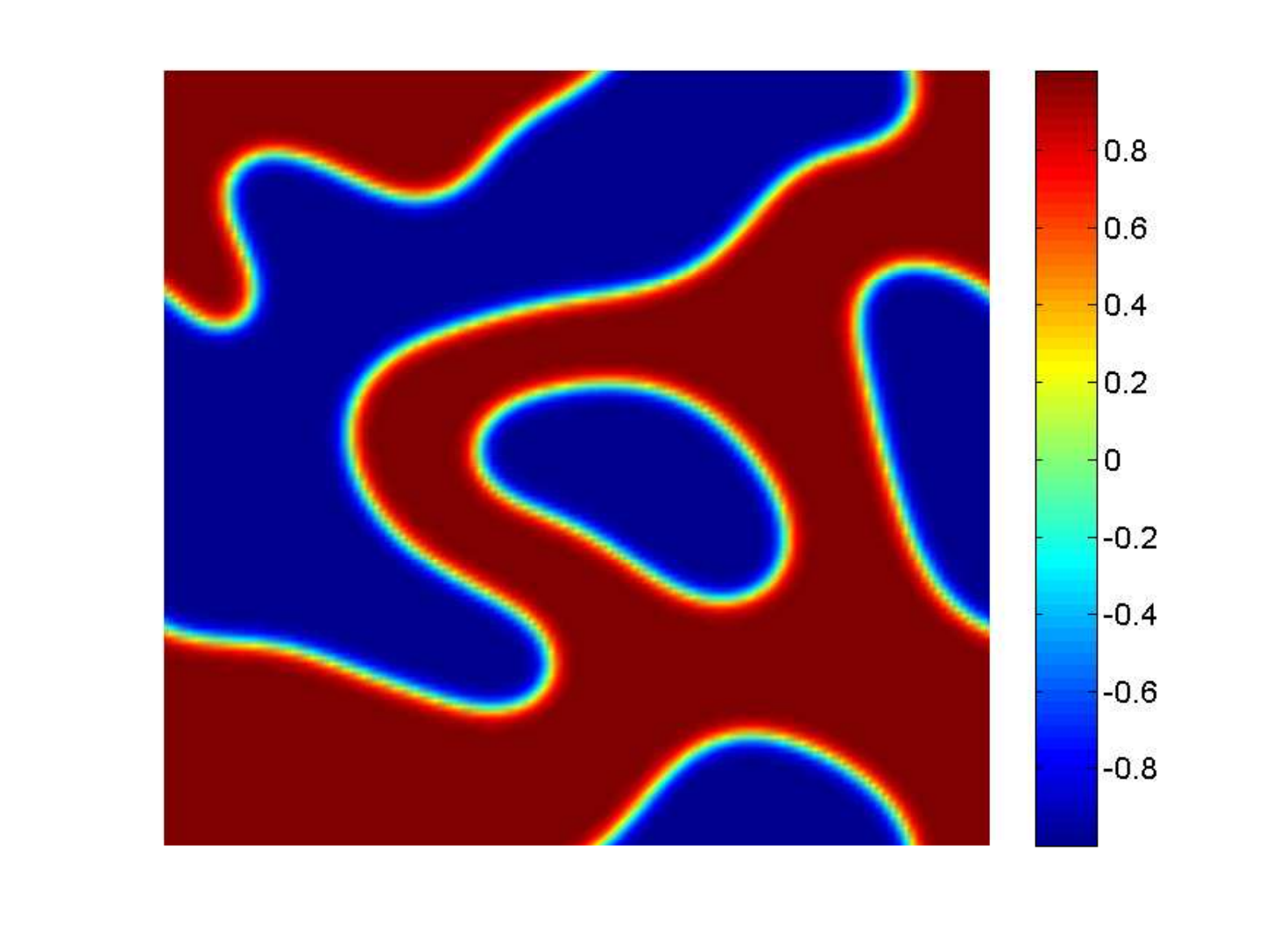}
\includegraphics[width=2.0in]{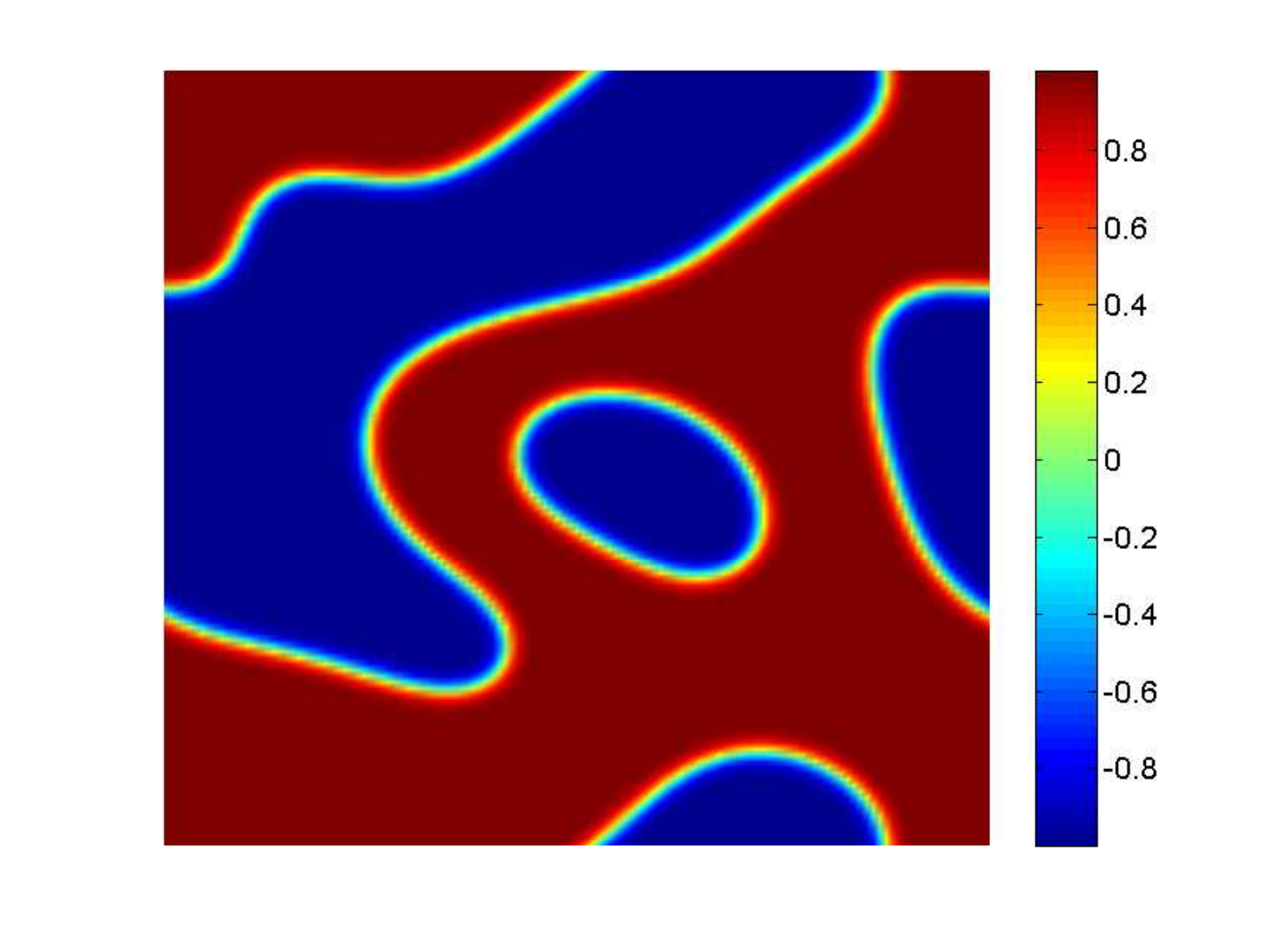}
\includegraphics[width=2.0in]{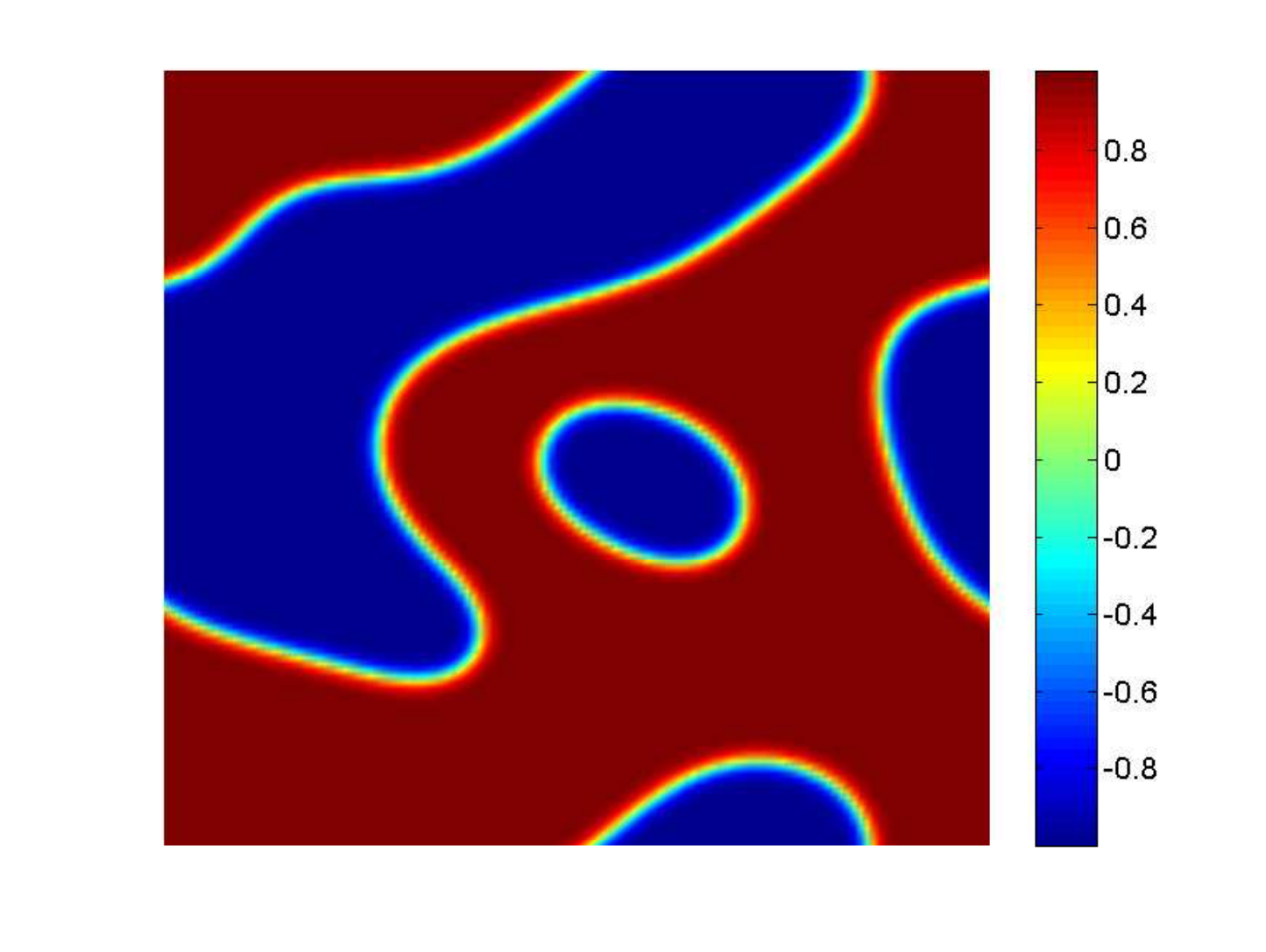}\\
\caption{Solution snapshots of coarsening dynamics
   of Allen-Cahn equation using adaptive time strategy at $t=1, 10, 20, 50,80,100$, respectively.}
\label{Snapshots-Dynamics}
\end{figure}

\begin{figure}[htb!]
\centering
\includegraphics[width=2.0in]{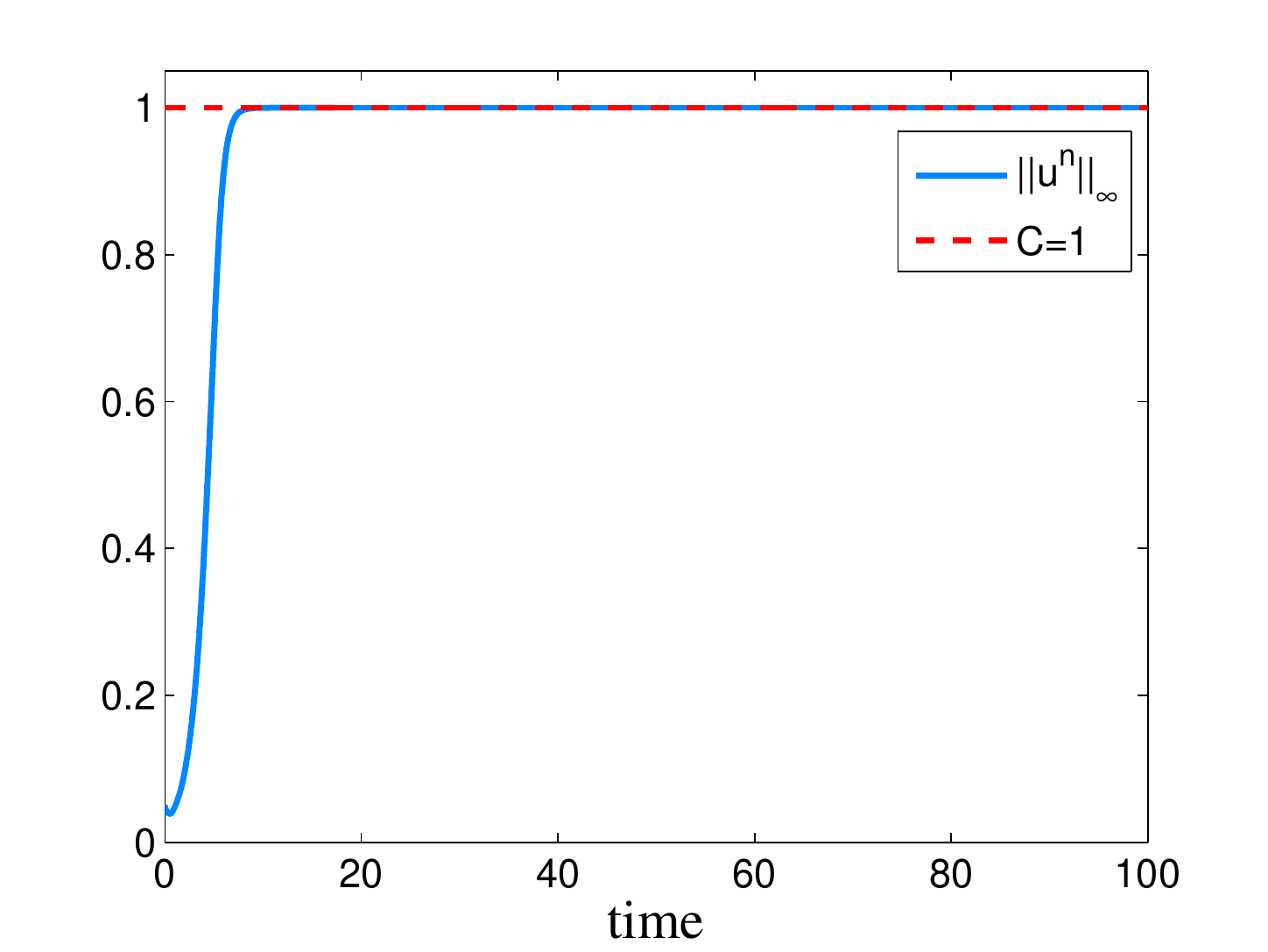}
\includegraphics[width=2.0in]{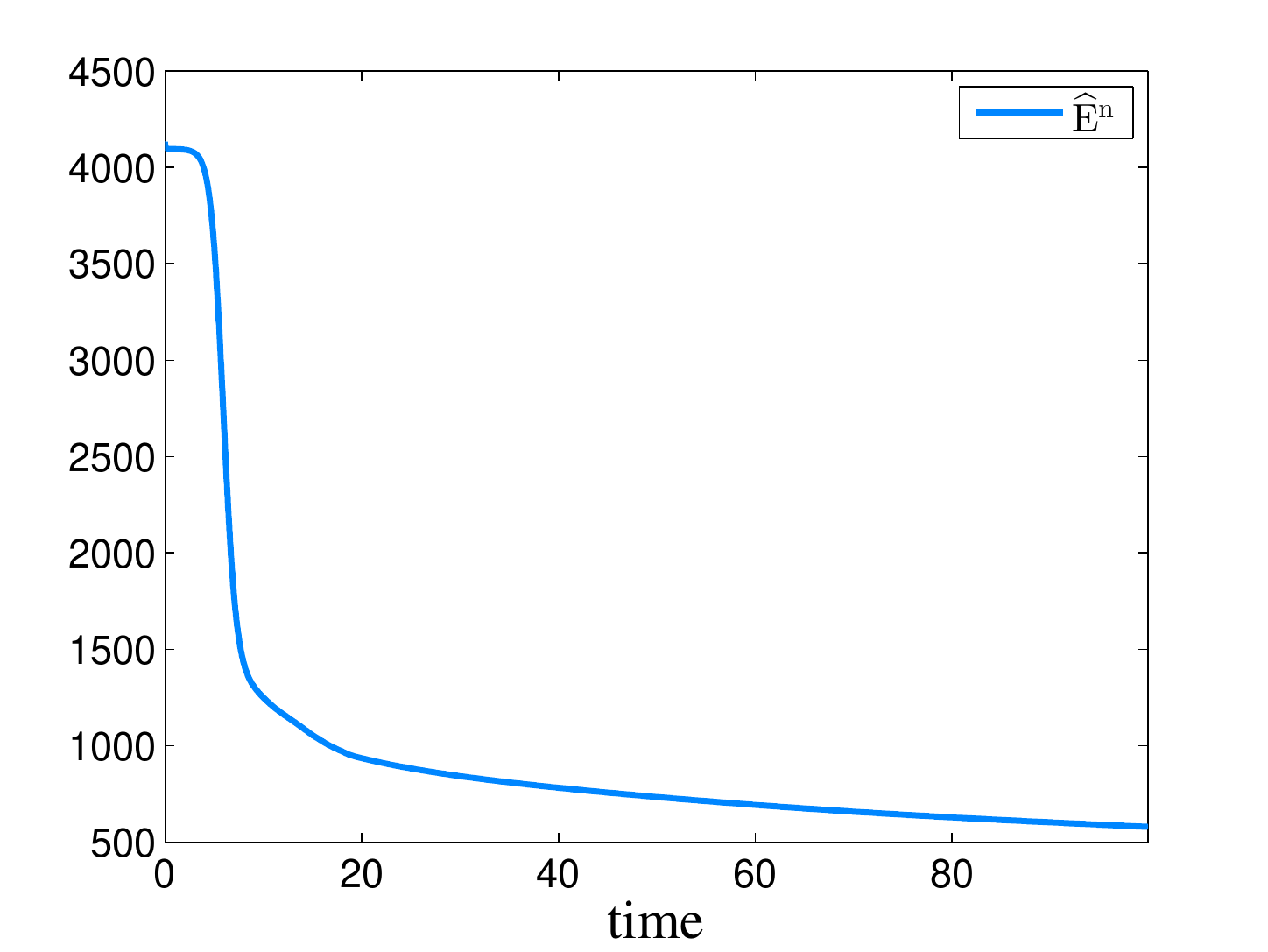}
\includegraphics[width=2.0in]{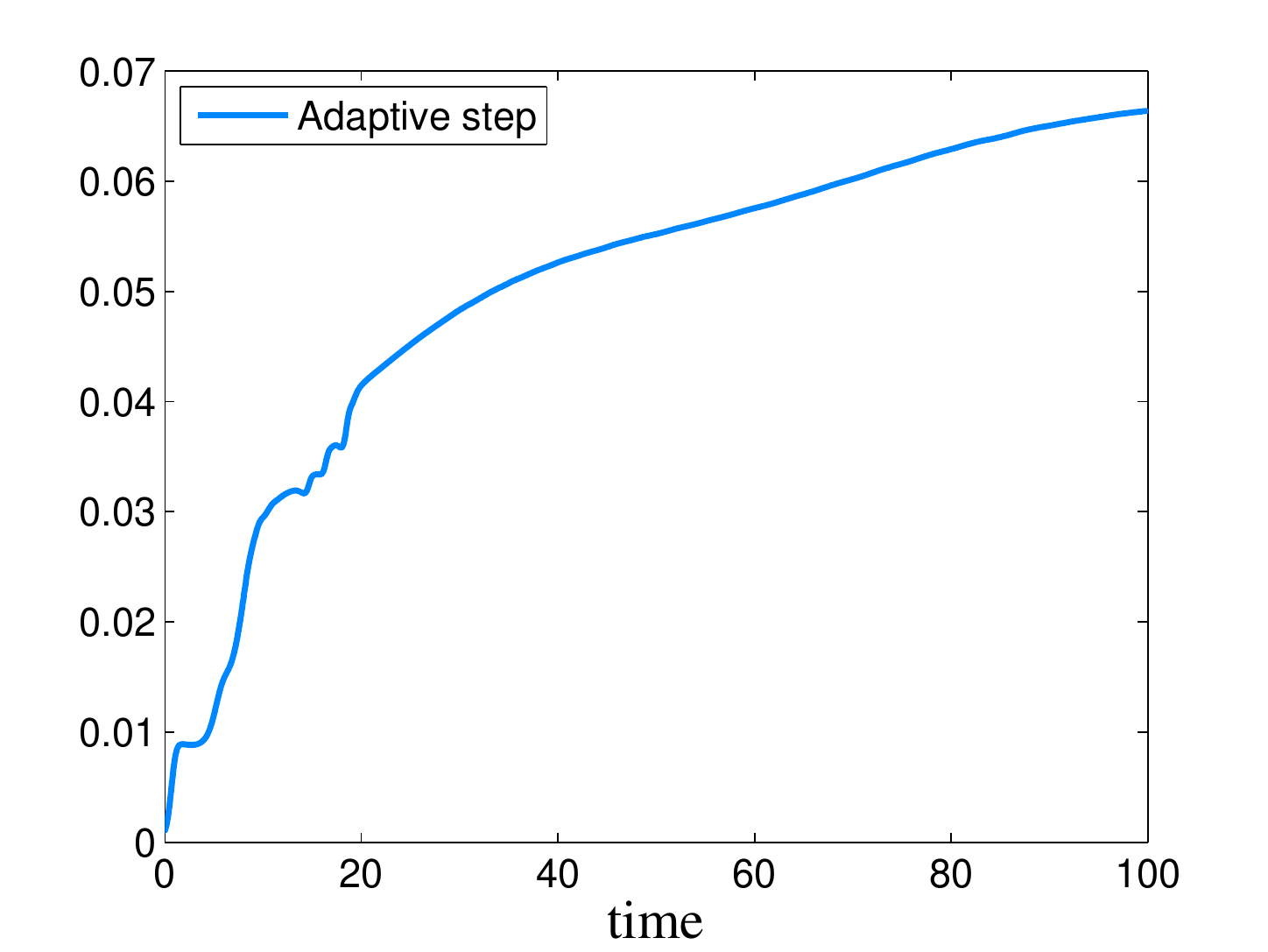}\\
\caption{Evolutions of Maximum norm (left), energy (middle) and adaptive time steps (right) of
  coarsening  dynamics of Allen-Cahn equation using adaptive time strategy.}
\label{Dynamics-MaxPrinciple-Energy}
\end{figure}

We next investigate the coarsening dynamic of the Allen-Cahn model by using adaptive BDF2 scheme
incorporated with the adaptive algorithm until $T=100$.
Figure \ref{Snapshots-Dynamics} shows the time evolution of the coarsening dynamic.
As can be seen at $t=1$, the microstructure is relatively fine
and contains a large number of grains.
As time evolves, the coarsening dynamic through migration of the phase boundaries,
decomposition and merging procedure can be observed.
Also, as a consequence the number of the grains becomes smaller with time.
The corresponding discrete maximum norm, energy and adaptive time step are plotted in Figure
\ref{Dynamics-MaxPrinciple-Energy}, where we observe that the maximum value of the numerical solutions are bounded by $1,$
the discrete energy decays monotonically, and the adaptive strategy is rather effective.

\section{Concluding remarks}
\setcounter{equation}{0}

This work is concerned with fully discretized numerical schemes for the Allen-Cahn equations. The main task of this work is to establish the energy stability, maximum principle and convergence analysis for the second-order BDF scheme with variable time steps. It is of practical importance to allow the use of variable time steps as the solutions of the Allan-Cahn equations may undergo different time regimes and require fine or coarse time steps accordingly. Of course, the ratio of the meshsize may increase or decrease smoothly in order to retain numerical stability. Consequently, some upper bounds may apply in the practical computations.

In this work, by using an appropriate energy method we have shown that the nonuniform BDF2 scheme preserve the energy dissipation law under a mild time ratio constraint. By using a kernel recombination and complementary technique, we show that the discrete maximum principle holds for the nonuniform BDF2 scheme under the time ratio constraint $r_k < 1 +\sqrt{2}$, which coincides with the Grigorie 's zero-stability condition . This maximum-principle preserving result seems very new for second-order time discretizations to the Allen-Cahn equation. This discrete maximum principle allows us to obtain the error estimates without any Lipschitz assumptions on the nonlinear bulk force. With the use of KRC technique and a new Gronwall inequality, the second-order rate of convergence in the maximum norm is finally established.

It is expected that the KRC technique developed in this work can be used to deal with more general nonlinear problems. One challenging topic is to develop nonuniform BFD2 type schemes for the time-fractional phase field equations \cite{LiaoTangZhou:2019,TangYuZhou:2019}. As the time-fractional operators require solution information at all time levels, the use of variable time steps seems more important in practice. On the technical side, it is of interests to see if the ratio constraints {\bf S1} and {\bf S0} are optimal or not.


\section*{Acknowledgements}
The authors would like to thank Dr. Bingquan Ji for his kind help on numerical simulations.


\end{document}